\documentclass{birkjour}
\numberwithin{equation}{section}
\newcommand{\CQFD}{\nolinebreak\hfill\rule{2mm}{2mm}\medbreak\par}
\newtheorem{theorem}{Theorem}[section]
\newtheorem{lemma}{Lemma}[section]
\theoremstyle{definition}
\newtheorem{definition}{Definition}[section]

\theoremstyle{remark}
\newtheorem{remark}{Remark}[section]

\begin{document}

%-------------------------------------------------------------------------
% editorial commands: to be inserted by the editorial office
%
%\firstpage{1} \volume{228} \Copyrightyear{2004} \DOI{003-0001}
%\seriesextraline{This is the Concrete Title of this Book\br H.E. R and S.T.C. W, Eds.}
%\firstpage{1}

\title[Order of the Central Moments of the Length of the LCS]
 {On the Order of the Central Moments of the Length of the Longest Common Subsequences in Random Words}

%----------Author 1
\author[Christian Houdr\'e]{Christian Houdr\'e}

\address{%
School of Mathematics, \\
Georgia Institute of Technology, \\
Atlanta, Georgia, 30332-0160}

\email{houdre@math.gatech.edu}

\thanks{This work was supported in part by the grant \#246283 from the Simons 
Foundation and by a Simons Foundation Fellowship, grant \#267336. Many thanks to 
the Laboratory MAS of the \'Ecole Centrale Paris and to the LPMA of the
Universit\'e Pierre et Marie Curie (Paris VI) for their hospitality and support while
part of this research was carried out.}
%----------Author 2
\author{Jinyong Ma}
\address{School of Mathematics, \\
Georgia Institute of Technology, \\
Atlanta, Georgia, 30332-0160}
\email{jinyma@gmail.com}
\subjclass{60K35; 60C05; 05A05}

\keywords{Longest common subsequence, $r$-th central moment, last passage percolation, 
Burkholder inequality, Efron-Stein inequality.}

\date{September 18, 2015}

\begin{abstract}
\noindent We investigate the order of the $r$-th, $1\le r < +\infty$, central moment
of the length of the longest common subsequence of two independent random words of size
$n$ whose letters are identically distributed and independently drawn from a finite alphabet.
When all but one of the letters are drawn with small probabilities, which depend on the size of the
alphabet, a lower bound is shown to be of order $n^{r/2}$.
This result complements a generic upper bound also of order $n^{r/2}$.
\end{abstract}

%%% ----------------------------------------------------------------------
\maketitle
%%% ----------------------------------------------------------------------
%\tableofcontents
\section{Introduction and statements of results}

Let $X=(X_i)_{ i\ge 1}$ and $Y=(Y_i)_{ i \ge 1}$ be two independent sequences of
iid random variables taking their values in a finite alphabet
$\mathcal{A}_m=\{\alpha_1,\alpha_2,\dots,\alpha_m\}$, $m\ge 2$, with $\mathbb{P}(X_1=\alpha_k)=\mathbb{P}(Y_1=\alpha_k)=p_k$, $k=1,2, \dots, m$.
Let also $LC_n$ be the length of the longest common subsequence of
the random words $X_1 \cdots X_n$ and $Y_1 \cdots Y_n$, {\it i.e.}, $LC_n:=LC_n(X_1 \cdots X_n;Y_1 \cdots Y_n)$ is the largest
$k$ such that there exist $1\le i_1< i_2 <\dots <i_k \le n$ and
$1\le j_1< j_2 <\dots <j_k \le n$, with
$X_{i_s}=Y_{j_s}$, $s=1,\dots,k$.

The study of the asymptotic behavior of $LC_n$ has a long history starting with the
well known result of Chv\'{a}tal and Sankoff \cite{CS} asserting that
\begin{equation}\label{lcmean}
\lim_{n \to \infty}\frac{\mathbb{E}LC_n}{n}=\gamma_m^*.
\end{equation}
However, to this day, the exact value of $\gamma_m^*$ (which depends on the distribution of
$X_1$ and on the size of the alphabet) is still unknown even
in "simple cases" such as for uniform Bernoulli random variables.
This first asymptotic result was sharpened by Alexander (\cite{A}) who showed that
\begin{equation}
\gamma_m^*n -K_{A}\sqrt{n \log n} \le \mathbb{E}LC_n \le \gamma_m^*n,
\end{equation}
where $K_{A}>0$ is a constant depending neither on $n$ nor on the distribution of $X_1$.
Next, Steele \cite{St} was the first to investigate the order of the variance
proving, in particular, that $VarLC_n\le n$. However, finding the order of the lower 
bound is more illusive.
For Bernoulli random variables and in various instances where there is a strong "bias"
such as high asymmetry or mixed common and increasing subsequence problems,
the lower bound is also shown to be of order $n$ (\cite{HLM}, \cite{HM}, \cite{LM}).
The uniform case is still unresolved and tight lower variance estimates
seem to be lacking (however, see \cite{AHM1}, \cite{AHM}, where a
situation "as close as we want" to uniformity is treated).

Below, starting with a generic upper bound, we investigate the order of the $r$-th, $r\ge1$, 
central moment
of $LC_n$ in case of finite alphabets (of course, as far as the order
is concerned only the case $1 \le r \le 2$ is really of interest for this lower bound).

The upper bound obtained in \cite{St} relies on an asymmetric version of
the Efron-Stein inequality which can be viewed as a
tensorization property of the variance.
The symmetric Efron-Stein inequality has seen a generalization, due to
Rhee and Talagrand \cite{RT}, to the
$r$-th moment where it is, in turn, viewed as a consequence of
Burkholder's square function inequality. As described next, in the asymmetric case,
a similar extension also holds thus providing
a generic upper bound on the $r$-th central moment of $LC_n$.
First, let $S:\mathbb{R}^n \to \mathbb{R}$ be a Borel function
and let $(Z_i)_{1\le i \le n}$ and $(\hat{Z}_i)_{1\le i \le n}$ be two independent families
of iid random variables having the same law.  Now, and with suboptimal
notation, let $S=S(Z_1,Z_2,\dots, Z_n)$,
and let $S_i=S(Z_1,Z_2,\dots, Z_{i-1}, \hat{Z}_i,Z_{i+1}, \dots,Z_n)$,
$1\le i \le n$.  Then, as shown next, for any $r \ge 2$,
\begin{equation}\label{mar-lcs}
\|S-\mathbb{E}S\|_r:=\left(\mathbb{E}|S-\mathbb{E}S|^r\right)^{{1}/{r}}\le \frac{r-1}{2^{1/r}}
\left(\sum_{i=1}^n \|S-S_i\|_r^2\right)^{1/2}.
\end{equation}

Indeed, for $i =1,\dots, n$, let $\mathcal{F}_i=\sigma(Z_1, \dots, Z_i)$
be the $\sigma$-field generated by $Z_1, \dots, Z_i$, let $\mathcal{F}_0=\{\Omega, \emptyset\}$
be trivial, and let $d_i:=\mathbb{E}(S|\mathcal{F}_i)-\mathbb{E}(S|\mathcal{F}_{i-1})$.
Thus, $(d_i,\mathcal{F}_i)_{1 \le i \le n}$ is a martingale differences sequence and from
Burkholder's square function inequality, with optimal constant, e.g., see \cite{OS}, for $r \ge 2$,
\begin{align}\label{mar1}
\|S-\mathbb{E}S\|_r =\left\|\sum_{i=1}^n d_i \right\|_r \!\le (r-1)\!\left\|\left(\sum_{i=1}^n
d_i^2\right)^{1/2}\right\|_r
\!\le (r-1)\!\left(\sum_{i=1}^n \|d_i^2\|_{r/2}\!\right)^{1/2}.
\end{align}
Moreover, and as in \cite{RT},
letting $\mathcal{G}_i=\sigma(Z_1, Z_2, \dots, Z_i, \hat{Z}_i)$,
$1\le i \le n$,
\begin{align}
\mathbb{E}|S-S_i|^r &=\mathbb{E}(\mathbb{E}(|S-S_i|^r|\mathcal{G}_i)) \notag\\
&\ge \mathbb{E}(|\mathbb{E}(S|\mathcal{G}_i)-\mathbb{E}(S|\mathcal{F}_{i-1})
+\mathbb{E}(S_i|\mathcal{F}_{i-1})-\mathbb{E}(S_i|\mathcal{G}_i)|^r)\notag\\
&:=\mathbb{E}|U+V|^r,\label{mar3}
\end{align}
where
$U=\mathbb{E}(S|\mathcal{G}_i)-\mathbb{E}(S|\mathcal{F}_{i-1})$
and $V=\mathbb{E}(S_i|\mathcal{F}_{i-1})-\mathbb{E}(S_i|\mathcal{G}_i)$.
But, given $\mathcal{F}_{i-1}$, $U$ and $V$ are independent, with moreover
$\mathbb{E}(U|\mathcal{F}_{i-1})=\mathbb{E}(V|\mathcal{F}_{i-1})=0$
and $\mathbb{E}|U|^r=\mathbb{E}|V|^r=\mathbb{E}|d_i|^r$, thus,
\begin{align}\label{mar4}
\mathbb{E}|U+V|^r=\mathbb{E}(\mathbb{E}(|U+V|^r|\mathcal{F}_{i-1}))
\ge \mathbb{E}|U|^r+\mathbb{E}|V|^r= 2\mathbb{E}|d_i|^r,
\end{align}
using the calculus inequality, valid for any $r\ge2$, $u \in \mathbb{R}$ and $v \in \mathbb{R}$, $|u+v|^r \ge |u|^r + r  sign(u)|u|^{r-1} v + |v|^r$, and taking conditional expectations. Combining \eqref{mar1}, \eqref{mar3} and \eqref{mar4} gives \eqref{mar-lcs}.

Next, apply \eqref{mar-lcs} to $LC_n$ viewed as a function of the $2n$ random variables $X_1,\dots, X_n,Y_1, \dots, Y_n$ and note, at first, that
replacing $X_i$ (resp. $Y_i$)
by an independent copy ${\hat X}_i$ (resp. ${\hat Y}_i$), changes $|LC_n-LC_n(X_1 \cdots \hat{X}_i \cdots X_n ; Y_1 \cdots Y_n)|$
(resp. $|LC_n-LC_n(X_1 \cdots X_n; Y_1 \cdots \hat{Y}_i \cdots Y_n)|$)
by at most $1$.  Thus, following Steele~\cite{St} and for each $i=1,\dots, n$,
\begin{align}
&||LC_n-LC_n(X_1 \cdots \hat{X}_i \cdots X_n ; Y_1 \cdots Y_n)||_r^2 \notag\\
& = \left(\mathbb{E}(|LC_n-LC_n(X_1 \cdots \hat{X}_i
\cdots X_n ; Y_1 \cdots Y_n)|^r \mathbf{1}_{X_i \ne \hat{X}_i})\right)^{2/r} \notag\\
& \le \left(\mathbb{P}(X_i \ne \hat{X}_i)\right)^{2/r}
= \left(1- \sum_{k=1}^m p_k^2\right)^{2/r}.\label{new-1}
\end{align}
Combining \eqref{new-1}, and its version for $(Y_i)_{1 \le i \le n}$,
with \eqref{mar-lcs} yields, for any $r\ge 2$,

\begin{equation}\label{E:theorem1-upperbound}
 \mathbb{E}|LC_n - \mathbb{E}LC_n|^r \le \frac{(r-1)^r}{2}
 \left(1- \sum_{k=1}^m p_k^2\right) (2n)^{{r}/{2}},
\end{equation}
which further yields,\label{E:theorem1-upperbound2}
$$\mathbb{E}|LC_n - \mathbb{E}LC_n|^r \le \left(\left(1- \sum_{k=1}^m p_k^2\right)n\right) ^{{r}/{2}},$$
for any $0< r \le 2$, by the Cauchy-Schwarz inequality.

Therefore, \eqref{E:theorem1-upperbound} provides an upper bound whose order could also be obtained, in a simpler way, by integrating out the tail inequality given via Hoeffding's exponential martingale inequality. Let us now state the main result of
the paper which provides a lower bound on the $r$-th central moment
of $LC_n$, when all but one of the
symbols are drawn with very small probabilities.

\begin{theorem}\label{T:theorem1}
Let $1\le r < +\infty$, and let $(X_i)_{ i \ge 1}$ and $(Y_i)_{ i \ge 1}$ be two
independent sequences of iid random variables with values
in $\mathcal{A}_m=\{\alpha_1,\alpha_2,\dots,\alpha_m\}$, $m\ge2$, such that $\mathbb{P}(X_1=\alpha_k)=p_k$, $k=1,2, \dots, m$. Further, let $j_0 \in \{1,\dots,m\}$ be such that $\max_{ j\ne j_0}p_j \le \min(2^{-2} e^{-5} K_m/m, K_m/2m^2)$, where $K_m=\min(K, 1/800m)$ and $K=2^{-4} 10^{-2} e^{-67}$. Then, there exists a constant $C>0$ depending on $r$, $m$, $p_{j_0}$
and $\max_{j\ne j_0}p_j$, such that, for all $n \ge 1$,
\begin{equation}\label{E:theorem1}
\mathbb{M}_r(LC_n):=\mathbb{E}\left|LC_n-\mathbb{E}LC_n\right|^{r} \ge C n^{\frac{r}{2}}.
\end{equation}
\end{theorem}
An estimate on the constant $C$ present in \eqref{E:theorem1} is given in Remark~\ref{LCS:remark1}.

In contrast to \cite{HM}, \cite{HLM} or \cite{LM} which deal only with binary words, our results
are proved for alphabets of arbitrary, but fixed size $m$,
and are thus novel in that context as well
even for the variance, i.e., $r=2$.
Moreover, our results are no longer existential, but provide precise constants depending on the alphabet size.
As well known, e.g., see \cite{AHM1}, \cite{AHM}, the LCS problem is a last
passage percolation (LPP) problem with strictly increasing paths and dependent weights
and, therefore, in our context, the order of the variance is linear.
For the LPP problem with independent weights the variance is conjectured to
be sublinear. In view of \eqref{E:theorem1-upperbound} and \eqref{E:theorem1}, it is
tempting to conjecture,
and we do so, that when properly centered (by $\gamma_m^* n$) and normalized (by $\sqrt n$),
asymptotically, $LC_n$ has a normal component. (The limiting law is in fact normal, see \cite{HU}.) This conjecture might
appear surprising since in LPP with independent weights different limiting laws
are conjectured and have been proved to be such in the closely related
Bernoulli matching model \cite{MN}. It should finally also be noted that, as seen
in \cite{BH}, with another closely related model, the order $n^{r/2}$ on the central
moments does not guarantee normal convergence, but nevertheless a normal component is present.

As for the content of the rest of paper, Section~\ref{section1} presents a proof of
Theorem~\ref{T:theorem1} which relies on a key preliminary result, Theorem~\ref{T:theorem2},
whose proof is given in Section~\ref{section2}.

\section{Proof of Theorem~\ref{T:theorem1}}\label{section1}

The strategy of proof to obtain the lower bound is to first represent $LC_n$ as a
random function of the number of most probable letters $\alpha_{j_0}$.
In turn, this random function locally satisfies a reversed Lipschitz condition which ultimately gives the lower bound in Theorem~\ref{T:theorem1}.
This methodology extends, modifies and simplifies (and at times corrects)
the binary strategy of proof of \cite{HM} or \cite{LM} providing also a more quantitative result.

To start, and as in \cite{HM}, pick a letter equiprobably at random from all the
non-$\alpha_{j_0}$ letters in either one of the two finite sequences, of length $n$,
$X$ or $Y$ ({\it Throughout the paper, by finite sequences $X$ and $Y$, of length $n$, it
is meant that $X=(X_i)_{1\le i \le n}$ and $Y=(Y_i)_{1 \le i \le n}$}).
Next, change it to the most probable letter $\alpha_{j_0}$ and call the two
new finite sequences $\tilde{X}$ and $\tilde{Y}$.
Then the length of the longest common subsequence
of $\tilde{X}$ and $\tilde{Y}$, denoted
by $\widetilde{LC_n}$, tends, on an event of high probability, to be larger
than $LC_n$.  This is the content of the following theorem which is proved in the
next section.

\begin{theorem}\label{T:theorem2}
Let the hypothesis of Theorem~\ref{T:theorem1} hold.
Then, for all $n \ge 1$, there exists a set $\mathcal{B}_n
\subset \mathcal{A}_m^n \times \mathcal{A}_m^n$,
such that,
\begin{equation}
\mathbb{P}\left((X,Y)\in \mathcal{B}_n\right)
\ge 1-125 \exp\left(-\frac{ n(\max_{ j\ne j_0}p_j)^6}{5}\right),
\end{equation}
and such that for all $(x,y) \in \mathcal{B}_n$,
\begin{equation}
\mathbb{P}(\widetilde{LC_n}-LC_n=1|X=x, Y=y)\ge \frac{K}{m},
\end{equation}
\begin{equation}
\mathbb{P}(\widetilde{LC_n}-LC_n=-1|X=x, Y=y)\le \frac{K}{2m},
\end{equation}
where $K=2^{-4} 10^{-2} e^{-67}$.
\end{theorem}

As already mentioned, the proof of Theorem~\ref{T:theorem2} is given in the next section, let us
nevertheless indicate how it leads to the
lower bound on $\mathbb{M}_r(LC_n)$ given in Theorem~\ref{T:theorem1}. In fact, the arguments leading to the conclusion of Theorem~\ref{T:theorem1} remain valid under any hypotheses for which the conclusions of Theorem~\ref{T:theorem2} remain valid.

{\it From now on, assume without loss of generality that
$p_1>1/2$ and that $p_2= \max_{2\le j\le m}p_j$, so that $\alpha_1$ is the most probable letter and $\alpha_2$ the second most probable one.}

To begin with, let us present a few definitions.
For the two finite random sequences $X=(X_i)_{1\le i \le n}$ and $Y=(Y_i)_{1 \le i \le n}$,
let $N_1$ be the total number of letters $\alpha_1$ present in both sequences, i.e.,
$N_1$ is a binomial random variable with parameters $2n$ and $p_1$.
Next, by induction, define a finite collection of pairs of finite random sequences
$(X^k,Y^k)_{0 \le k \le 2n}$, which are independent of $X$ and $Y$, and therefore independent of $N_1$, as follows: First, let $X^0=(X^0_i)_{1 \le i \le n}$
and $Y^0=(Y^0_i)_{1 \le i \le n}$ be independent, with $X^0_i$ and $Y^0_i$, $i=1,\dots, n$,
iid random variables with values in $\{\alpha_2, \dots, \alpha_m\}$ and such that
$\mathbb{P}(X_1^0=\alpha_k)=\mathbb{P}(Y_1^0=\alpha_k)=p_k/(1-p_1)$, $2\le k \le m$.
In other words, $X^0$ and $Y^0$ are two independent finite sequences of iid random variables
whose joint law is the law of $((X,Y)|N_1=0)$.  Once ($X^k, Y^k$) is defined,
let ($X^{k+1}, Y^{k+1}$) be the pair of finite random sequences obtained by taking
(pathwise) with equal probability, one letter from all the
letters $\alpha_2, \alpha_3,\dots, \alpha_m$ in the pair $(X^k , Y^k)$ and replacing it
with $\alpha_1$, and for this path iterating the process till $k=2n$.
Clearly, for $1\le k \le 2n-1$, $X^k$ and $Y^k$ are not independent,
while $(X_i^{2n}, Y_i^{2n})_{1 \le i \le n}$ is a deterministic sequence made up only of the letter $\alpha_1$.

Rigorously, the random variables can be defined as follows: let $\Omega$ be our underlying space, and  let $\Omega^{2n+1}$ be its $(2n+1)$-fold Cartesian product. For each
$\omega = (\omega_0, \omega_1, \dots, \omega_{2n}) \in \Omega^{2n+1}$ and
$0 \le k \le 2n$, $(X^k(\omega), Y^k(\omega))$ only depends on
$\omega_0, \omega_1, \dots, \omega_k$.  Then, $(X^{k+1}(\omega), Y^{k+1}(\omega))$ is obtained
from $(X^k(\omega), Y^k(\omega))$ by replacing with equal probability
any non-$\alpha_1$ letter by $\alpha_1$,
while the choice of the non-$\alpha_1$ letter to be replaced in $(X^k(\omega) , Y^k(\omega))$ is
determined by $\omega_{k+1}$.

Next, let $LC_n(k)$ denote the length of the longest common subsequence
of $X^k$ and $Y^k$ (with a slight abuse of notation and
terminology with the identification of finite sequences and words).
The lemma below shows that $(X^k,Y^k)$ has the same law as
$(X,Y)$ conditional on $N_1=k$, and therefore the law of $LC_n(k)$ is the same as
the conditional law of $LC_n$ given $N_1=k$.

\begin{lemma}\label{E:lemma01}
For any $k=0, 1, \dots, 2n$,
\begin{equation}\label{E:lemma1}
(X^k, Y^k)\overset{d}{=}((X,Y)|N_1=k),
\end{equation}
and moreover,
\begin{equation}\label{E:lemma12}
(X^{N_1}, Y^{N_1})\overset{d}{=}(X,Y),
\end{equation}
where $\overset{d}{=}$ denotes equality in distribution.
\end{lemma}
$\textbf{Proof.}$ The proof is by induction on $k$. By definition,
$(X^0, Y^0)$ has the same law as $(X,Y)$ conditional on $N_1=0$.
For any $(\alpha_{j_1},\dots,\alpha_{j_{2n}}) \in \mathcal{A}_m^n \times \mathcal{A}_m^n$, let
$$q_{\ell}=\left|\left\{1\le i \le 2n: \ \alpha_{j_i}=\alpha_{\ell}\right\}\right|,$$
$1 \le \ell \le m$.
Now assume that \eqref{E:lemma1} is true for $k$, {\it i.e.}, assume that for any
$(\alpha_{j_1},\dots, \alpha_{j_{2n}}) \in \mathcal{A}_m^n \times \mathcal{A}_m^n$, with $q_1=k$,
\begin{equation}
\mathbb{P}\left((X_1^k, \dots, X_n^k, Y_1^k, \dots, Y_n^k)=(\alpha_{j_1},\dots,\alpha_{j_{2n}})\right)=\binom{2n}{k}^{-1}\prod_{\ell=2}^{m}\left(\frac{p_\ell}{1-p_1}\right)^{q_{\ell}}.
\end{equation}
Then, for any $(\alpha_{j_1},\dots,\alpha_{j_{2n}}) \in \mathcal{A}_m^n \times \mathcal{A}_m^n$,
with $q_1=k+1$,
\begin{multline}\label{LCS:lemma1}
\mathbb{P}\left((X_1^{k+1}, \dots, X_n^{k+1}, Y_1^{k+1}, \dots, Y_n^{k+1})=(\alpha_{j_1},
\dots,\alpha_{j_{2n}})\right)=\\
\sum_{i=1}^{k+1}\!\mathbb{P}\!\left(\!(X_1^{k+1}, \dots, X_n^{k+1}, Y_1^{k+1},
\dots, Y_n^{k+1})=(\alpha_{j_1},\dots,\alpha_{j_{2n}})|B_i^{k+1}\right)\!\mathbb{P}(B_i^{k+1}),
\end{multline}
where $B_i^{k+1}$, $1\le i \le k+1$, is the event that the $i$-th $\alpha_1$ in
$(\alpha_{j_1},\dots,\alpha_{j_{2n}})$ is changed from a non-$\alpha_1$ letter when passing
from $(X^k,Y^k)$ to $(X^{k+1},Y^{k+1})$. (Conditional on $B_i^{k+1}$, the $i$-th $\alpha_1$ in
$(\alpha_{j_1},\dots,\alpha_{j_{2n}})$ could have been changed from any letter in
$\{\alpha_2,\alpha_3,\dots,\alpha_m\}$.)  Assuming this
$\alpha_1$ has been changed, say, from $\alpha_s$, $2 \le s \le m$, the corresponding
probability is given by:
$$\mathbb{P}\left((X^k,Y^k)=(\alpha_{j_1},\dots,\alpha_s,\dots,\alpha_{j_{2n}})\right)=\binom{2n}{k}^{-1}  \prod_{\ell=2}^{m}\left(\frac{p_\ell}{1-p_1}\right)^{q_{\ell}}\left(\frac{p_s}{1-p_1}\right),$$
where, above, $\alpha_s$ takes the place of the $i$-th $\alpha_1$ in the sequence
$(\alpha_{j_1},\dots,\alpha_{j_{2n}})$. Thus,
\begin{align}
&\mathbb{P}\left((X_1^{k+1}, \dots, X_n^{k+1}, Y_1^{k+1}, \dots, Y_n^{k+1})=(\alpha_{j_1},\dots,\alpha_{j_{2n}})|B_i^{k+1}\right)\mathbb{P}(B_i^{k+1})\notag\\
&=\binom{2n}{k}^{-1}\prod_{\ell=2}^{m}\left(\frac{p_\ell}{1-p_1}\right)^{q_{\ell}}
\left(\sum_{s=2}^m \frac{p_{s}}{1-p_1}\right) \frac{1}{2n-k},\notag
\end{align}
which when incorporated into \eqref{LCS:lemma1}, gives
\begin{multline}
\mathbb{P}\left((X_1^{k+1}, \dots, X_n^{k+1}, Y_1^{k+1}, \dots, Y_n^{k+1})=(\alpha_{j_1},
\dots,\alpha_{j_{2n}})\right)\\
=\binom{2n}{k+1}^{-1}\prod_{\ell=2}^{m}\left(\frac{p_\ell}{1-p_1}\right)^{q_{\ell}},
\end{multline}
finishing the proof of the first part of the lemma.

Next, from \eqref{E:lemma1} and the independence of $N_1$ and $\{(X^k, Y^k)\}_{0\le k\le 2n}$,
for any $(u,v) \in \mathbb{R}^{n}\times \mathbb{R}^{n}$,
\begin{align}
\mathbb{E}\left(e^{i<u,X>+i<v,Y>}\right)&=\sum_{k=0}^{2n}\mathbb{E}
\left(e^{i<u,X>+i<v,Y>}|N_1=k\right)\mathbb{P}\left(N_1=k\right)\notag\\
&=\sum_{k=0}^{2n}\mathbb{E}\left(e^{i<u,X^k>+i<v,Y^k>}\right)\mathbb{P}\left(N_1=k\right)\notag\\
&=\sum_{k=0}^{2n}\mathbb{E}\left(e^{i<u,X^{k}>+i<v,Y^{k}>}|N_1=k\right)\mathbb{P}\left(N_1=k\right)\notag\\
&=\sum_{k=0}^{2n}\mathbb{E}\left(e^{i<u,X^{N_1}>+i<v,Y^{N_1}>}|N_1=k\right)\mathbb{P}\left(N_1=k\right)\notag\\
&=\mathbb{E}\left(e^{i<u,X^{N_1}>+i<v,Y^{N_1}>}\right),\notag
\end{align}
finishing the proof of the lemma.\CQFD

Let now $LC_n(N_1)$ be the length of the longest common subsequence of $X^{N_1}$ and $Y^{N_1}$.
The above lemma implies that $LC_n$ and $LC_n(N_1)$ have the same law and, therefore,
\begin{equation}\label{equal}
\mathbb{M}_r(LC_n(N_1))=\mathbb{M}_r(LC_n).
\end{equation}
To lower bound the right hand side of \eqref{equal} (and to
prove Theorem~\ref{T:theorem1}) the following simple inequality will prove useful.
\begin{lemma}\label{L:lemma2.2}
Let $f:Dom\rightarrow\mathbb{Z}$ satisfy a local reversed Lipschitz condition, i.e., let $\ell \ge 0$
and let $f$ be such that for any $i,j \in D$ with $ j\ge i+ \ell$,
$$f(j)-f(i)\ge c(j-i),$$
for some $c>0$.  Let $T$ be a $Dom$-valued random variable
with $\mathbb{E}|f(T)|^r<+\infty$, $r \ge 1$, then
\begin{equation}
\mathbb{M}_r(f(T))\ge  \left(\frac{c}{2}\right)^r \left(\mathbb{M}_r(T)-\ell^r\right).
\end{equation}
\end{lemma}
$\textbf{Proof.}$ Let $r\ge 1$, and let $\widehat{T}$ be an independent copy of $T$.
First, and clearly,
$\mathbb{M}_r(T)\le \mathbb{E}(|T-\widehat{T}|^r)\le 2^r\mathbb{M}_r(T)$.
Hence,
\begin{align}
\mathbb{M}_r(f(T))& \ge \frac{1}{2^r} \mathbb{E}(|f(T)-f(\widehat{T})|^r)\notag\\
& \ge \left(\frac{c}{2}\right)^r \left(\mathbb{E}(T-\widehat{T})^r\mathbf{1}_{T-\widehat{T}\ge \ell}
+\mathbb{E}(\widehat{T}-T)^r\mathbf{1}_{\widehat{T}-T\ge \ell}\right)\notag\\
& \ge  \left(\frac{c}{2}\right)^r \left(\mathbb{E}|T-\widehat{T}|^r-\ell^r\right)\notag\\
& \ge  \left(\frac{c}{2}\right)^r \left(\mathbb{M}_r(T)-\ell^r\right).\notag
\end{align} \CQFD

The above lemma will prove useful in providing a lower
bound on $\mathbb{M}_r(\!LC_n(N_1))$ by showing that, after removing
the randomness of $LC_n(\cdot)$, $LC_n(\cdot)$ satisfies a local
reversed Lipschitz condition. To do so, for a random variable $U$ with finite $r$-th moment and
for a random vector $V$, let $\mathbb{M}_r(U|V):=\mathbb{E}
\left(|U-\mathbb{E}\left(U|V\right)|^r\big| V\right)$.
Clearly, by convexity and the conditional Jensen's inequality,
\begin{align}\label{E:inequality1}
 \mathbb{M}_r(U|V) & \le 2^r \left(\mathbb{E}\left(|U-\mathbb{E}U|^r \big| V \right)/2 + \mathbb{E}\left(|\mathbb{E}(U|V)-\mathbb{E}U|^r \big| V \right)/2 \right)\notag\\
 & \le 2^r \mathbb{E}\left(|U-\mathbb{E}U|^r \big| V \right)
\end{align}
and so, for any $n \ge 1$,
\begin{align}
\mathbb{M}_r(LC_n(N_1))&\ge \frac{1}{2^r}
\mathbb{E}(\mathbb{M}_r(LC_n(N_1)|(LC_n(k))_{0\le k \le 2n}))\notag\\
&= \frac{1}{2^r} \int_{\Omega} \mathbb{M}_r(LC_n(N_1)|(LC_n(k))_{0\le k \le 2n}(\omega))
\mathbb{P} (d \omega)\notag\\
& \ge \frac{1}{2^r} \int_{O_n} \mathbb{M}_r(LC_n(N_1)|(LC_n(k))_{0\le k \le 2n}(\omega))
\mathbb{P} (d \omega)\label{E:lowerbound1},
\end{align}
where for each $n\ge 1$,
\begin{equation}
O_n:= \bigcap_{\substack{i,j \in I\\j\ge i+\ell(n) }}\left\{LC_n(j)-LC_n(i)\ge  \frac{K}{4m}(j-i)\right\},
\end{equation}
where $K$ is given in Theorem~\ref{T:theorem2} and where $\ell(n) \ge 0$ is to be chosen later.
(Of course, above and everywhere, intersections, unions and sums are taken
over countable sets of integers.) In words, on the event $O_n$ the random function $LC_n$ has a
slope of at least $K/4m$, when restricted to the interval $I$
and when $i$ and $j$ are at least $\ell(n)$ apart from each other.

Since $N_1$ is independent of $(LC_n(k))_{0 \le k \le 2n}$, and
from \eqref{E:inequality1}, for each $\omega \in \Omega$,
\begin{align}\label{E:lowerbound2}
&\mathbb{M}_r(LC_n(N_1)|(LC_n(k))_{0\le k \le 2n}(\omega))\notag\\
& \ge \frac{1}{2^r}\mathbb{M}_r(LC_n(N_1)|(\!LC_n(k))_{0\le k \le 2n}(\omega),
\mathbf{1}_{N_1\in I}=1)\mathbb{P}(N_1 \in I | (LC_n(k))_{0\le k \le 2n}(\omega))\notag\\
& = \frac{1}{2^r}\mathbb{M}_r(LC_n(N_1)|(LC_n(k))_{0\le k \le 2n}(\omega),
\mathbf{1}_{N_1\in I}=1)\mathbb{P}(N_1 \in I),
\end{align}
where
\begin{equation}\label{E:interval}
I=\left[2np_1 -\sqrt{2n(1-p_1)p_1 },2np_1 +\sqrt{2n(1-p_1)p_1 }\right].
\end{equation}
Again, for each $\omega \in O_n$, from Lemma \ref{L:lemma2.2}, and since $N_1$ is independent of $(LC_n(k))_{0 \le k \le 2n}$,
\begin{multline}\label{E:lowerbound5-new}
 \mathbb{M}_r(LC_n(N_1)|(LC_n(k))_{0\le k \le 2n} (\omega), \mathbf{1}_{N_1\in I }=1)\\
 \ge   \left(\frac{K}{8m}\right)^r \left(\mathbb{M}_r(N_1| \mathbf{1}_{N_1 \in I }=1)-\ell(n)^r\right).
\end{multline}

Now, \eqref{E:lowerbound1}, \eqref{E:lowerbound2} and \eqref{E:lowerbound5-new} lead to
\begin{equation}\label{E:lowerbound4-new}
\mathbb{M}_r(LC_n(N_1))
\ge \frac{1}{4^r}\!\left(\frac{K}{8m}\right)^r\!\!\left(\mathbb{M}_r(N_1|
\mathbf{1}_{N_1 \in I }=1)-\ell(n)^r\right)\mathbb{P}(N_1 \in I) \mathbb{P}(O_n),
\end{equation}
and it remains to estimate each one of the three terms on the right hand side of
\eqref{E:lowerbound4-new}. By the Berry-Ess\'een inequality, and all $n\ge 1$,
\begin{equation}\label{LCS:BE}
\left|\mathbb{P}(N_1 \in I)- \frac{1}{\sqrt{2 \pi}}
\int_{-1}^1 e^{-\frac{x^2}{2}} dx\right|\le \frac{1}{\sqrt{2np_1(1-p_1)}}.
\end{equation}
Moreover,
\begin{align}\label{condition1}
& \mathbb{M}_r(N_1|\mathbf{1}_{N_1 \in I }=1)\notag\\
&= \mathbb{E}(|N_1-2np_1+2np_1-\mathbb{E}(N_1|\mathbf{1}_{N_1 \in I }=1)|^r | \mathbf{1}_{N_1 \in I }=1)\notag\\
&\ge \left|\mathbb{E}(|N_1-2np_1|^r | \mathbf{1}_{N_1 \in I }=1)^{1/r}-|2np_1-\mathbb{E}(N_1|\mathbf{1}_{N_1 \in I }=1)|\right|^r,
\end{align}
and
\begin{align}
&|\mathbb{E}(N_1|\mathbf{1}_{N_1 \in I }=1)-2np_1| \notag\\
&=\sqrt{2np_1(1-p_1)}\left|\mathbb{E}\left(\frac{N_1-2np_1}
{\sqrt{2np_1(1-p_1)}}\Big| \mathbf{1}_{N_1 \in I }=1 \right)\right|\notag\\
& = \sqrt{2np_1(1-p_1)}\frac{\left|\! F_n(1)-\Phi(1)+F_n(-1)-\Phi(-1)-\int_{-1}^{1}(F_n(x)-\Phi(x))dx\right|}{\mathbb{P}(N_1 \in I)}\notag\\
& \le \sqrt{2np_1(1-p_1)} \frac{4 \max_{x \in [-1,1]}|F_n(x)-\Phi(x)|}
{\mathbb{P}(N_1 \in I)}\notag\\
& \le  \frac{2}{\int_{-1}^1 e^{-\frac{x^2}{2}}dx/{\sqrt{2 \pi}} - 1/{\sqrt{2np_1(1-p_1)}}},\label{condition2}
\end{align}
where $F_n$ is the distribution functions of ${(N_1-2np_1)}/{\sqrt{2np_1(1-p_1)}}$, while
$\Phi$ is the standard normal one. Likewise,
\begin{align}
&\mathbb{E}(|N_1-2np_1|^r | \mathbf{1}_{N_1 \in I }=1) \notag\\
&\ge (2n p_1(1-p_1))^{r/2}\frac{\int_{-1}^1 |x|^r d\Phi(x) -4\max_{x \in [-1,1]}|F_n(x)-\Phi(x)|}{\mathbb{P}(N_1 \in I)}   \notag\\
& \ge (2np_1(1-p_1))^{r/2} \frac{\int_{-1}^1 |x|^r e^{-\frac{x^2}{2}}dx-2\sqrt{\pi}/\sqrt{np_1(1-p_1)}}{\int_{-1}^1 e^{-\frac{x^2}{2}} dx+\sqrt{\pi}/\sqrt{np_1(1-p_1)}}.\label{condition3}
\end{align}
Next, \eqref{condition1}-\eqref{condition3} lead to:
\begin{align}
&\mathbb{M}_r(N_1|\mathbf{1}_{N_1 \in I }=1) \notag\\
&\ge \Bigg| (2np_1(1-p_1))^{\frac{1}{2}}\left( \frac{\int_{-1}^1 |x|^r e^{-\frac{x^2}{2}}dx-2\sqrt{\pi}/\sqrt{np_1(1-p_1)}}{\int_{-1}^1 e^{-\frac{x^2}{2}} dx+
\sqrt{\pi}/\sqrt{np_1(1-p_1)}}\right)^{\frac{1}{r}}\notag\\
&\quad \quad \quad \quad\quad\quad\quad\quad\quad\quad\quad\quad\quad -\frac{2}{\int_{-1}^1
e^{-\frac{x^2}{2}}dx/{\sqrt{2 \pi}} - 1/{\sqrt{2np_1(1-p_1)}}}\Bigg|^r.\label{E:lowerbound7}
\end{align}
Finally, assuming Theorem~\ref{T:theorem2}, the
estimates \eqref{E:lowerbound4-new}-\eqref{E:lowerbound7} combined with
the estimate on $\mathbb{P}(O_n)$ obtained in the next lemma give
the lower bound \eqref{E:theorem1}, whenever $33m^2\log n /K^2 \le \ell(n)
\le K_1 \sqrt{n}$ (where $K_1$ is given and estimated in Remark~\ref{LCS:remark1}).

\begin{lemma}
For $m\ge 2$, let $K_m=\min(K, 1/800m)$ where $K=2^{-4} 10^{-2} e^{-67}$, and let $p_2\le \min(2^{-2} e^{-5} K_m/m, K_m/2m^2)$.  Then, for all $n \ge 1$,
\begin{equation}\label{LCS:inlemma2.3}
\mathbb{P}(O_n) \ge 1-\left(500\sqrt{\pi}e^2 n  \exp\left(-\frac{ np_2^6}{5}\right)+2n \exp\left(-\frac{K^2\ell(n)}{32m^2}\right)\right).
\end{equation}
\end{lemma}
$\textbf{Proof.}$ Let $A_n:=\{(X,Y)\in \mathcal{B}_n\}$ and let $A_n^k:=\{(X^k,Y^k)\in \mathcal{B}_n\}$. Then,
\begin{equation}\label{E:leitem2}
\mathbb{P}\left(\left(\bigcap_{k \in I} A_n^k\right)^c\right)
\le \sum_{k \in I}\mathbb{P}\left(\left(A_n^{k}\right)^c\right)= \sum_{k \in I}\mathbb{P}
\left(A_n^c|N_1=k\right)\le \sum_{k \in I}\frac{\mathbb{P}(A_n^c)}{\mathbb{P}(N_1=k)},
\end{equation}
by Lemma~\ref{E:lemma01}.
Next, by Stirling's formula in the form,
$$\sqrt{2\pi} n^{n+\frac{1}{2}} e^{-n+\frac{1}{12n+1}}< n! < \sqrt{2 \pi} n^{n+\frac{1}{2}} e^{-n + \frac{1}{12n}},$$
for all $k \in I$ and $n\ge 1$,
\begin{align}
\mathbb{P}(N_1=k)&=\binom{2n}{k} p_1^k (1-p_1)^{2n-k}\notag\\
&\ge \frac{1}{\sqrt{2 \pi} e^2}
\frac{(2n)^{2n+1/2}}{k^{k+1/2}(2n-k)^{2n-k+1/2}} p_1^k (1-p_1)^{2n-k}\notag\\
&:= \gamma (k,n,p_1).\notag
\end{align}
Hence, for all $k \in I$ and $p_1 \ge 3/4$ (which holds true since $p_2 \le K/m$), from the property of the probability mass function of the binomial distribution,
\begin{align}
&\mathbb{P}(N_1=k) \notag\\
& \ge \min \left(\mathbb{P}(N_1=2np_1 -\lfloor \sqrt{2n(1-p_1)p_1 } \rfloor), \mathbb{P}(N_1=2np_1 +\lfloor\sqrt{2n(1-p_1)p_1 }\rfloor)\right)\notag\\
&\ge \min\left(\gamma \left(2np_1 -\lfloor \sqrt{2n(1-p_1)p_1 }\rfloor,n,p_1\right), \gamma \left(2np_1 +\lfloor \sqrt{2n(1-p_1)p_1}\rfloor,n,p_1\right)\right)\notag\\
& \ge \frac{1}{2\sqrt{2 \pi} e^2\sqrt{n}}.
\end{align}
This last inequality in conjunction with \eqref{E:leitem2} and Theorem~\ref{T:theorem2},
gives
\begin{equation}
\mathbb{P}\left(\left(\bigcap_{k \in I} A_n^k\right)^c\right)\le 4\sqrt{\pi}e^2 n \mathbb{P}(A_n^c)\le 500\sqrt{\pi}e^2 n  \exp\left(-\frac{ np_2^6}{5}\right).
\end{equation}
Next, for each $n\ge 1$, letting
\begin{equation}\label{E:defofdelta}
\Delta_{k+1}=
\begin{cases}
LC_n(k+1)-LC_n(k), &\text{when $A_n^k$ holds,}\\
1, &\text{otherwise,}\\
\end{cases}
\end{equation}
it follows from Theorem~\ref{T:theorem2} that,
\begin{equation}\label{E:leitem1}
\mathbb{E}(\Delta_{k+1}|X^k,Y^k)\ge \frac{K}{2m}.
\end{equation}
Now, for each $k=0,1,\dots, 2n$, let $\mathcal{F}_k:=\sigma(X^0,Y^0, \dots, X^k,Y^k)$,
be the $\sigma$-field generated by $X^0,Y^0, \dots, X^k,Y^k$.  Clearly,
$(\Delta_{k}-\mathbb{E}(\Delta_{k}|\mathcal{F}_{k-1}), \mathcal{F}_k)_{1\le k \le 2n}$
forms a martingale differences sequence and since $-1\le \Delta_k\le 1$, Hoeffding's martingale inequality gives, for any $i<j$,
\begin{equation}
\mathbb{P}\left(\sum_{k=i+1}^j \left(\Delta_{k}-\mathbb{E}(\Delta_{k}|\mathcal{F}_{k-1})\right)
< -\frac{K}{4m}(j-i)\right)\le \exp\left(-\frac{K^2(j-i)}{32m^2}\right).
\end{equation}
Moreover, from \eqref{E:leitem1}, $\sum_{k=i+1}^j\mathbb{E}(\Delta_{k}|X^{k-1},Y^{k-1})
\ge {K}(j-i)/2m$, and therefore
\begin{align}\label{E:leitem3}
\mathbb{P}\left(\sum_{k=i+1}^j \Delta_{k} \le \frac{K}{4m}(j-i)\right)& \le
\mathbb{P}\left(\sum_{k=i+1}^j \left(\Delta_{k}-\mathbb{E}(\Delta_{k}|\mathcal{F}_{k-1})\right) <
-\frac{K}{4m}(j-i)\right)\notag\\
&\le \exp\left(-\frac{K^2(j-i)}{32m^2}\right).
\end{align}
For each $n \ge 1$, let now
$$O_n^{\Delta}=\bigcap_{\substack{i,j \in I\\j \ge i+\ell(n)}}\left\{\sum_{i+1}^j \Delta_k
\ge \frac{K}{4m}(j-i)\right\},$$
then, from \eqref{E:leitem3}
\begin{equation}
\mathbb{P}\left(\left(O_n^{\Delta}\right)^c\right)\le \sum_{\substack{i,j \in I\\j \ge i+\ell(n)}}
\mathbb{P}\left(\sum_{i+1}^j \Delta_k < \frac{K}{4m}(j-i)\right)
\le 2n \exp\left(-\frac{K^2\ell(n)}{32m^2}\right).
\end{equation}
From the very definition of $\Delta_k$ in \eqref{E:defofdelta},
$\bigcap_{k \in I} A_n^k \cap O_n^{\Delta} \subset O_n,$ and therefore
\begin{align}
\mathbb{P}\left(\left(O_n\right)^c\right)&\le
\mathbb{P}\left(\left(\bigcap_{k \in I} A_n^k\right)^c\right)+\mathbb{P}
\left(\left(O_n^{\Delta}\right)^c\right)\notag\\
&\le 500\sqrt{\pi}e^2 n  \exp\left(-\frac{ np_2^6}{5}\right)
+2n \exp\left(-\frac{K^2\ell(n)}{32m^2}\right).
\end{align}
 \CQFD

\begin{remark}\label{LCS:remark1}
The reader might wonder how to estimate the constant $C$ in
Theorem~\ref{T:theorem1}. In view of \eqref{equal}, the right hand side of \eqref{E:lowerbound4-new} needs to be lower bounded. Letting $n \ge p_2^{-12}+m^8$, together with \eqref{LCS:BE}, \eqref{E:lowerbound7} and \eqref{LCS:inlemma2.3} yield to:
$$\mathbb{P}(N_1 \in I)\ge \frac{1}{2},\quad \mathbb{P}(O_n) \ge \frac{1}{2},$$
and
$$\mathbb{M}_r(N_1|\mathbf{1}_{N_1 \in I }=1) \ge e^{-\frac{1}{2}} 2^{-(1+r)}(1+r)^{-1}(n(1-p_1))^{\frac{r}{2}}.$$
Moreover, choosing
$$\ell(n)=2^{(-1-r-\frac{1}{r})}e^{-\frac{1}{2r}}(n(1-p_1))^{\frac{1}{2}}\left(\frac{1}{1+r}\right)^{\frac{1}{r}}:= K_1 \sqrt{n},$$
in \eqref{E:lowerbound4-new}, gives:
$$\mathbb{M}_r(LC_n) \ge 2^{-4-6r}(1+r)^{-1}e^{-1/2}K^r m^{-r}(1-p_1)^{r/2}n^{r/2}.$$
Letting $C_{1}=2^{-4-6r}(1+r)^{-1}e^{-1/2}K^r m^{-r}(1-p_1)^{r/2}$, and
$$C_2=\min_{n \le p_2^{-12}+m^8}\frac{\mathbb{M}_r(LC_n) }{n^{r/2}}
\le \frac{(r-1)^r}{2} 2^{{r}/{2}} \left(1- \sum_{k=1}^m p_k^2\right),$$
by \eqref{E:theorem1-upperbound}, then one can
choose $C=\min(C_1, C_2)$ in Theorem~\ref{T:theorem1}.
\end{remark}

\section{Proof of Theorem \ref{T:theorem2}}\label{section2}

\subsection{Description of alignments}

Let us begin with an example. Let $\mathcal{A}_3=\{\alpha_1,\alpha_2,\alpha_3\}$, with $\alpha_i = i$,  $i=1,2,3,$ and, say that
\begin{equation}\label{sequences}
X=121313111211,\quad
Y=111311112112.
\end{equation}
An optimal alignment of $X$ and $Y$, i.e., an alignment corresponding to a LCS, is
\begin{equation}\label{A:1}
\begin{tabular}{l r r r r r r r r r r r r r}
\hline
1 & 2 &\ & 1 & 3 & 1 & 3     & 1 & 1 &1 &2 &1 &1 &\ \\
1 & \ &1 & 1 & 3 & 1 & \    & 1 & 1 &1 &2 &1 &1 &2 \\
\hline
\end{tabular}
\end{equation}
and another possible optimal alignment is
\begin{equation}\label{A:2}
\begin{tabular}{l r r r r r r r r r r r r r}
\hline
1 & 2 &1 & \ & 3 & 1 & 3 & 1 & 1 &1 &2 &1 &1 &\ \\
1 & \ &1 & 1 & 3 & 1 & \  & 1 & 1 &1 &2 &1 &1 &2 \\
\hline
\end{tabular}
\end{equation}
both corresponding to the LCS $1131111211$.

Comparing these two optimal alignments, it is clear that the way the letters $\alpha_1$ are
aligned, between the aligned non-$\alpha_1$ letters, is not important as
long as a maximal number of such letters $\alpha_1$ are aligned.
Therefore, in general, it is enough to describe which non-$\alpha_1$ letters are aligned and to
assume that between pairs of aligned non-$\alpha_1$ letters a maximal number of letters
$\alpha_1$ are aligned. In other words, we can identify the
two optimal alignments \eqref{A:1} and \eqref{A:2}
as the same.

Next, let a {\it cell}, be either the beginning of an alignment till and including, if any,
its first pair of aligned non-$\alpha_1$ letter, or be a part of an alignment between pairs
of aligned non-$\alpha_1$ letters.

For example, the alignment \eqref{A:1} can be decomposed
into two cells $C(1)$ and $C(2)$ as
\begin{equation}
\overbrace{
\begin{tabular}{l r r r r}
\hline
1 & 2 &\ & 1 & 3\\
1 & \ &1 & 1 & 3\\
\hline
\end{tabular}}^{C(1),\ v_1=-1}
\overbrace{
\begin{tabular}{l r r r r r }
\hline
1 & 3  & 1 & 1 &1 &2\\
1 & \  & 1 & 1 &1 &2\\
\hline
\end{tabular}}^{C(2),\ v_2=0}
\begin{tabular}{l r r}
\hline
1 & 1 & \  \\
1 & 1 & 2  \\
\hline
\end{tabular}
\end{equation}
where, moreover, each $v_i$ denotes the difference between the number of
letters $\alpha_1$ in the $X$-strand and the $Y$-strand of
the cell $C(i)$, $i=1,2$. For the alignment \eqref{A:1},
this gives the representation $v=(v_1,v_2)=(-1,0)$.
Another optimal alignment is via $v=(v_1,v_2)=(0,-1)$ corresponding
to another LCS, namely 1113111211:
\begin{equation}
\overbrace{
\begin{tabular}{l r r r r r}
\hline
1 & 2  & 1 & 3  & 1 & 3\\
1 & \  & 1 & \  & 1 & 3\\
\hline
\end{tabular}}^{C(1),\ v_1=0}
\overbrace{
\begin{tabular}{l r r r r   }
\hline
1 & 1 & 1 &\ &2\\
1 & 1 & 1 &1 &2\\
\hline
\end{tabular}}^{C(2),\ v_2=-1}
\begin{tabular}{l r r}
\hline
1 & 1 & \  \\
1 & 1 & 2  \\
\hline
\end{tabular}
\end{equation}
Note that any alignment has a cell-decomposition with a
corresponding finite vector of differences.
(With the convention that when no non-$\alpha_1$ letters are aligned, then the
alignment has no cell.)

Let $X=X_1X_2\cdots X_n$ and $Y=Y_1Y_2\cdots Y_n$ be given.
As just conveyed, any alignment has a cell-decomposition with an associated
vector representation
$v:=(v_1,\dots, v_k)$ indicating the number of cells ($k$, here)
and the differences between the number of letters $\alpha_1$ in the $X$-strand and
the corresponding number in the $Y$-strand of each cell. Conversely,
any $v \in \mathbb{Z}^k$ corresponds to a, possibly empty, family of cell-decompositions.

Let us now turn to {\it optimality}.  First, clearly any optimal alignment
is made of, say, $k$ cells (recall also our convention above),
where within each cell a maximum number of letters $\alpha_1$ are aligned
and, if any, the optimal alignment also has a tail part (the part after the last cell, i.e.,
the part after the last aligned non-$\alpha_1$ letters)
where as many letter $\alpha_1$ as possible are aligned.
Therefore, such an optimal alignment is given via a
unique $v \in \mathbb{Z}^k$.  On the other hand,
every $v=(v_1,\dots,v_k)\in \mathbb{Z}^k$ also corresponds to a (possibly empty)
family of optimal alignments.
All of these optimal alignments have the same number of pairs of aligned
non-$\alpha_1$ letters where within each cell a maximal
number of letters $\alpha_1$ are aligned, and where moreover
as many letters $\alpha_1$
as possible are aligned after the pair of aligned non-$\alpha_1$-letters.
These optimal alignments corresponding to the same $v$
can differ in the way
the letters $\alpha_1$ are aligned within each cell and in the tail part.
It can also happen, and in contrast to the binary case, that one can align
different pairs of non-$\alpha_1$ letters, which can only happen when
no letters $\alpha_1$ are present between these different pairs of non-$\alpha_1$
letters. (Take, for example, $X=1321$ and $Y=2311$, then the optimal alignments
corresponding to $v\in \mathbb{Z}$ can align either the letter $2$ or the letter $3$.)
But in both cases such optimal alignments based on the same $v$ give
the same length for the corresponding longest common subsequences.
Therefore, we can identify all the optimal alignments in the family
associated with $v$ as a single one.
In other words, we identify each vector $v$ with an optimal alignment,
provided one exists, and vice-versa.

Writing $|v|$ for the number of coordinates of $v$, {\it i.e.}, $|v|=k$, if $v \in \mathbb{Z}^k$,
the cell-decomposition $\pi-\nu$
associated with $v=(v_1,\dots,v_k)\in \mathbb{Z}^k$ can now precisely be defined:

\begin{definition}\label{D:alignment}
Let $k\in \mathbb{N}, k \ge 1$ and let $v=(v_1,\dots, v_k)\in \mathbb{Z}^k$.
Let $\pi_{v}(0)=\nu_{v}(0)=0$, and for each $i=1, \dots, k$, let $(\pi_{v}(i),\nu_{v}(i))$
be any one of the smallest pair of integers $(s,t)$ (where $(s_1,t_1)\le(s_2,t_2)$ indicates
that $s_1\le s_2$ and $t_1\le t_2$) satisfying the following three conditions:
\begin{enumerate}
\item $\pi_{v}(i-1)<s$ and $\nu_{v}(i-1)<t$;
\item $X_s=Y_t\in\{\alpha_2,\dots, \alpha_m\}$;
\item the difference between the number of letters $\alpha_1$ in the
integer intervals $[\pi_{v}(i-1),s]$
and $[\nu_{v}(i-1),t]$ is equal to $v_{i}$.
\end{enumerate}
If for some $i= 1, \dots, k$, no such $(s,t)$ exists, then set $\pi_{v}(i)=\dots=\pi_{v}(k)=\infty$
and $\nu_{v}(i)=\dots=\nu_{v}(k)=\infty$.
\end{definition}

In other words, above, $\pi_{v}(i), \nu_{v}(i)$, $i=1, \dots, k$, are the indices corresponding
to the $i$-th aligned non-$\alpha_1$ pair in $v$.  For $i=1, \dots, k$,
the $i$-th cell, $C_v(i)$ is the
pair
$$C_v(i):=\left(X_{\pi_{v}(i-1)+1} \dots X_{\pi_{v}(i)};
Y_{\nu_{v}(i-1)+1} \cdots Y_{\nu_{v}(i)}\right),$$
and the cell $C_v(i)$ \emph{is called a $v_i$-cell.}

Let us further comment on the above definition, we actually defined a
greedy algorithm for each cell (each cell must be {\it minimal} meaning that
the cell ends as soon as all three conditions in Definition~\ref{D:alignment} are met).
For any optimal alignment, let us compare its cells with our minimal cells alignment.
If any, respectively denote the first two different cells by
$c_{i}^{opt}$ and $c_{i}^{min}$, $1\le i \le k$, since these cells
correspond to the same $v_i \in \mathbb{Z}$, they only differ in the number of pairs of
aligned letters $\alpha_1$.
From the definition of minimality, $c_{i}^{opt}$ contains more pairs of aligned
letters $\alpha_1$ than $c_{i}^{min}$.   These pairs of letters $\alpha_1$,
being of same number on the $X$-strand and $Y$-strand, can thus be pushed to next cell.
By iterating this push-procedure till the tail, then any optimal alignment can be
transformed into a minimal (optimal) alignment without reducing the length of
the common subsequence.
Thus an optimal alignment can always be transformed into a minimal (optimal) alignment.

With the above definition, we can let the alignment
associated to $v$ be any alignment (provided one exists) satisfying the following three
conditions:
\begin{enumerate}
\item $X_{\pi_{v}(i)}$ is aligned with $Y_{\nu_{v}(i)}$, for every $i=1,2,\dots,k$;
\item the number of aligned letters $\alpha_1$ in the cell $C_v(i)$, denoted by $S_v(i)$, is the minimum
number of letters $\alpha_1$ present in either
\!$X_{\pi_{v}(i-1)+1}\! \cdots\! X_{\pi_{v}(i)}$\! or $Y_{\nu_{v}(i-1)+1}\! \cdots\! Y_{\nu_{v}(i)}$;
\item after having aligned $X_{\pi_{v}(k)}$ with $Y_{\nu_{v}(k)}$, then align as many letters
$\alpha_1$ as possible and denote that number by $r_v$.
\end{enumerate}

From these definitions, for any $v \in \mathbb{Z}^k$, and if there exists a 
minimal cell-decomposition corresponding to $v$ exists,
then $\pi_v(k)\le n$ and $\nu_v(k)\le n$.  Such a $v$ is then
said to be \emph{admissible}.  Let $V$ denote the set of all admissible cell-decompositions, that is,
\begin{equation}\label{admissible}
V:=\left\{v\in \bigcup_{k = 1}^{\infty}\mathbb{Z}^k: \pi_{v}(|v|)\le n, \nu_{v}(|v|)\le n\right\}.
\end{equation}
Then, for every $v \in V$, and further for $|v|=0$ in case of no cell,
the length of the common subsequence corresponding to this
alignment is:
\begin{equation}
\Lambda C_v=|v|+\sum_{i=1}^{|v|}S_v(i)+r_v.
\end{equation}
Therefore the length of the longest common subsequence of $X$ and $Y$ can be expressed as:
\begin{equation}
LC_n=\max_{v \in V} \Lambda C_v,
\end{equation}
and, moreover, an alignment associated to an admissible $v$ is optimal
if and only if $\Lambda C_v=LC_n$.

\subsection{The effect of changing a non-$\alpha_1$ letter into $\alpha_1$}\label{sectiondef}

Again, the main idea behind Theorem~\ref{T:theorem2} is that, by changing a randomly picked non-$\alpha_1$
letter into $\alpha_1$, the length of the longest common subsequence is more likely to increase
by one than to decrease by one.  More precisely, conditional on the event
$A_n=\{(X,Y)\in \mathcal{B}_n\}$, the probability of an increase of $LC_n$ is at least
$K/m$ while the probability of a decrease is at most $K/2m$.
Let us illustrate this fact with another example. Let $X$ and $Y$ be given by,
\begin{equation}
X=112113112131, \ Y=131111111131,
\end{equation}
with optimal alignment:
\begin{equation}
\overbrace{
\begin{tabular}{l r r r r r r r r r r r r r}
\hline
1 & \ &1 & 2 & 1 &1 &3 &1 &1 &2 &1 &\ &\ &3 \\
1 & 3 &1 & \ & 1 &1 &\ &1 &1 &\ &1 &1 &1 &3 \\
\hline
\end{tabular}}^{C(1),\ v_1=-2}
\begin{tabular}{l }
\hline
1 \\
1 \\
\hline
\end{tabular}
\end{equation}
Above, there are 6 non-$\alpha_1$ letters, $X_3, X_6, X_9, X_{11}, Y_2, Y_{11}$, and each one has probability  $1/6$ to be picked and replaced by $\alpha_1$. Next, $X_3, X_6, X_9$ and $Y_2$ are not aligned with other letters but rather with gaps.
Moreover, since $X_3, X_6, X_9$ are on the top strand which contains a lesser number of
letters $\alpha_1$, picking one of them and replacing it leads to an increase of one in
the length of the LCS. On the other hand, since $X_{11}$ and $Y_{11}$ are aligned in this optimal
alignment, picking one of them and replacing it could potentially (but not necessarily) decrease the length of the LCS
by one. Finally, picking $Y_2$ may only potentially increase the length of the LCS by modifying
the alignment. In conclusion, in this example, by switching a randomly chosen non-$\alpha_1$ letter
into $\alpha_1$, the probability of an increase of the length of the LCS is at least $1/2$,
while the probability of a decrease is at most $1/3$.

To prove Theorem~\ref{T:theorem2}, we just need to prove that typically there exists
an optimal alignment such that:
\begin{enumerate}
\item Among all the non-$\alpha_1$ letters in $X$ and $Y$, the proportion which are on the
cell-strand with the smaller number of letters $\alpha_1$ is at least $K/m$.
\item Among all the non-$\alpha_1$ letters in $X$ and $Y$, the proportion which is
aligned is at most $K/2m$.
\end{enumerate}

Formally, let $v=(v_1, \dots, v_k) \in \mathbb{Z}^k$ be admissible.
For each $1 \le i \le k$, if $v_i\ne 0$, let $N_{v}^{-}(i)$ be the number of
non-$\alpha_1$ letters on the cell-strand of $C_v(i)$ with the lesser number of
letters $\alpha_1$, i.e., let
\begin{equation}
N_v^-(i)=
\begin{cases}
\sum_{j=\pi_v(i-1)+1}^{\pi_v(i)-1} \mathbf{1}_{X_j \in \{\alpha_2,\dots, \alpha_m\}} , &\text{if $v_i<0$,}\\
\sum_{j=\nu_v(i-1)+1}^{\nu_v(i)-1} \mathbf{1}_{Y_j \in \{\alpha_2,\dots, \alpha_m\}} , &\text{if $v_i>0$,}
\end{cases}
\end{equation}
while if $v_i= 0$, let $N_{v}^{-}(i)=0$.  Then, the total number of
non-$\alpha_1$ letters present on the cell-strands with the smaller number of 
letters $\alpha_1$ is equal to
\begin{equation}\label{weakside}
N_v^-:=\sum_{i=1}^{|v|} N_v^-(i).
\end{equation}
Let $N_i$ be the number of letters $\alpha_i$ in the two finite sequences
$X$ and $Y$, and let
\begin{equation}\label{totalnon1}
N_{>1}=\sum_{i=2}^m N_i.
\end{equation}
Next, let
\begin{multline}\notag
\mathcal{B}_n:=\left\{(x,y) \in \mathcal{A}_m^n \times \mathcal{A}_m^n:
\text{ there exists an optimal alignment } \text{ of }(x,y)  \right.\\
\left. \text{ with } |v| \ge 1,
n_v^-\ge K n_{>1}/m \text{ and } 2|v|\le K  n_{>1}/2m \right\},
\end{multline}
where, above, $n_v^-$ is the value of $N_v^-$ corresponding to $v$ and similarly for $n_{>1}$.
Clearly, $\mathcal{B}_n$ depends on $K$ and $m$.
Letting $A_n=\{(X,Y)\in \mathcal{B}_n\}$, our goal is now to prove that for some $\tilde{K}>0$,
independent of $n$, $\mathbb{P}\left(A_n\right)\ge 1-e^{-\tilde{K} n}$.

To continue, we need an optimal alignment having enough non-$\alpha_1$ letters in the
cell-strands with the smaller number of letters $\alpha_1$. However, for many optimal alignments,
most cells are zero-cells,  \emph{i.e.}, cells with the same number of letters $\alpha_1$ on both strands.
To bypass this hurdle, on an optimal alignment where most cells are zero-cells, some of the zero-cells are broken
up in order to create enough nonzero-cells while at the same time, maintaining the
optimality of the alignment after this breaking procedure.
Let us present this breaking operation on an example.  Take the two sequences
$$X=112113113,\quad \text{and} \  Y=112131113.$$
One of their optimal alignments is
\begin{equation}
\overbrace{
\begin{tabular}{l r r }
\hline
1 & 1 &2 \\
1 & 1 &2\\
\hline
\end{tabular}}^{C(1),\ v_1=0}
\overbrace{
\begin{tabular}{l r r r r r r}
\hline
1 &\ &1 & 3  & 1 & 1 &3 \\
1 &3 &1 & \  & 1 & 1 &3 \\
\hline
\end{tabular}}^{C(2),\ v_2=0}
\end{equation}
where both cells $C(1)$ and $C(2)$ are zero-cells. Now in the cell $C(2)$, $X_6$ and $Y_5$ are only
one position away from being aligned. Thus aligning them, instead of the pair $X_5$
and $Y_6$, breaks the cell $C(2)$ into two new cells $\tilde{C}(2)$ and
$\tilde{C}(3)$, with $\tilde{v}_2=1$ and $\tilde{v}_3=-1$.
The new optimal alignment is then:
\begin{equation}
\overbrace{
\begin{tabular}{l r r }
\hline
1 & 1 &2 \\
1 & 1 &2\\
\hline
\end{tabular}}^{\tilde{C}(1),\ \tilde{v}_1=0}
\overbrace{
\begin{tabular}{l r r }
\hline
1 &1 &3 \\
1 &\ &3 \\
\hline
\end{tabular}}^{\tilde{C}(2),\ \tilde{v}_2=1}
\overbrace{
\begin{tabular}{l r r r}
\hline
1 &1 &\ &3  \\
1 &1 &1 &3 \\
\hline
\end{tabular}}^{\tilde{C}(3),\ \tilde{v}_3=-1}
\end{equation}

The advantage of breaking up a zero-cell is that the resulting newly formed cells 
have different numbers of
letters $\alpha_1$ on each strand, thus $N_v^-$ tends to increase in this process
while the length of the common subsequence remains the same.
After applying this procedure and getting enough cells with different numbers of
letters $\alpha_1$ on the two strands, there is a high probability of finding
enough non-$\alpha_1$ letters on the strand with the smaller number
of letters $\alpha_1$.

The previous example leads to our next definition.
\begin{definition}
Let $k \in \mathbb{N}, k\ge 1$,
let $v \in \mathbb{Z}^k \cap V$, and for $i=1, \dots, k$, let $C_v(i)$ be any cell with
$v_i=0$. Then, $C_v(i)$ is said to be breakable if there exist $j$ and
$j^\prime$ such that:
\begin{enumerate}
\item $X_j=Y_{j^\prime}\in \{\alpha_2,\dots,\alpha_m\}$;
\item $\pi_v(i-1)<j<\pi_v(i)$ and $\nu_v(i-1)<j^\prime<\nu_v(i)$;
\item the difference between the number of letters $\alpha_1$ in
$$X_{\pi_v(i-1)+1}X_{\pi_v(i-1)+2}\cdots X_{j-1}
\text{  and  }Y_{\nu_v(i-1)+1}Y_{\nu_v(i-1)+2}\cdots Y_{j^\prime-1}$$
is plus or minus one.
\end{enumerate}
\end{definition}

\subsection{Probabilistic developments}
After the combinatorial developments of the previous sections, let us now bring forward some
probabilistic tools.
We start by introducing a useful way of constructing alignments
corresponding to a
given vector $v=(v_1, \dots, v_k)\in \mathbb{R}^k$.

For $1 \le i \le n$ and $2\le j \le m$, let $R_i^j$ (resp. $S_i^j$) be the number of
letters $\alpha_j$ between the $(i-1)$-th and $i$-th $\alpha_1$ in the infinite
sequence $(X_i)_{i\ge 1}$
(resp. $(Y_i)_{i\ge 1}$), with, of course, $R_1^j$ (resp. $S_1^j$) being the
number of letters $\alpha_j$ before the first $\alpha_1$.

Recall also, from Definition~\ref{D:alignment}, that in order to construct a
zero-cell, we use the
random time $T_0$, given by
\begin{equation}\label{stoppingtime}
T_0=\min_{2\le j \le m}T_0^j,
\end{equation}
where $T_0^j:=\min\{i=1,2,\ldots:\quad R_i^j\ne 0,S_i^j\ne 0\}$.
For a $-u$-cell ($u>0$), the random time is
\begin{equation}
T_{-u}=\min_{2\le j \le m}T_{-u}^j,
\end{equation}
where $T_{-u}^j:=\min\{i=1,2,\ldots:R_i^j\ne 0,S_{i+u}^j\ne 0\}$, and for a $u$-cell ($u>0$),
\begin{equation}
T_{u}=\min_{2\le j \le m}T_u^j,
\end{equation}
where $T_{u}^j:=\min\{i=1,2,\ldots:R_{i+u}^j\ne 0,S_{i}^j\ne 0\}$.
In other words, a cell with $v_i=u$ can be constructed in the following way:
Begin by keeping the first $u$ letters $\alpha_1$ in the $X$-strand, then align consecutive pairs
of letter $\alpha_1$ until meeting the first pair of the same non-$\alpha_1$ letter. 
(As previously argued, here different choices of pairs of the same non-$\alpha_1$ 
letter are possible, {\it i.e.}, if there are no letters $\alpha_1$ between 
different minimal pairs, but any pair will do if there is more than one choice.)

Let us find the law of $R_i^j$ and, to do so, let $R_i^{>1}=\sum_{j=2}^m R_i^j$
be the total number of non-$\alpha_1$ letters between the $(i-1)$-th and the
$i$-th $\alpha_1$.  Then, $R_i^{>1}+1$ is a geometric random variable with
parameter $p_1$, i.e., $\mathbb{P}(R_i^{>1}=k)=(1-p_1)^k p_1$,
$k=0,1,2,\dots$.  Moreover, conditionally on $R_i^{>1}$, $(R_i^{j})_{j=2}^m$ has a
multinomial distribution and therefore
\begin{align}
\mathbb{P}(R_i^j=k)&=\sum_{\ell=k}^{\infty}\mathbb{P}(R_i^j=k|R_i^{>1}
=\ell)\mathbb{P}(R_i^{>1}=\ell)\notag\\
&=\sum_{\ell=k}^{\infty}
\binom{\ell}{k}\left(\frac{p_j}{1-p_1}\right)^k \left(\frac{1-p_1-p_j}{1-p_1}\right)^{\ell-k}
(1-p_1)^{\ell}p_1 \notag\\
&=\left(\frac{p_1}{p_1+p_j}\right)\left(\frac{p_j}{p_1+p_j}\right)^k,
\end{align}
for $k=0,1,2,\dots$.  Thus, $R_i^j+1$ has a geometric distribution
with parameter ${p_1}/{(p_1+p_j)}$, $2 \le j \le m$.

To continue our probabilistic analysis, let us provide a rough lower bound on the
length of the LCS. First, aligning as many letters $\alpha_1$ as possible in $X$ and $Y$, would get
approximately a common subsequence of length $np_1 $, then aligning as many
letters $\alpha_2$ as possible without disturbing the already aligned
$\alpha_1$, would give an additional $\sum_{i=1}^{np_1 }\min\{R_i^2, S_i^2 \}$ aligned $\alpha_2$.
Moreover, since $R_i^2$ and $S_i^2$ are independent geometric random variables,
$\min\{R_i^2, S_i^2 \}+1$ is a geometric random variable with parameter
$1-(p_2/{(p_1+p_2)})^2$.  So, on average, the aligned letters
$\alpha_2$ contribute to the length of the LCS by an amount of:
$$np_1 \frac{p_2^2}{p_1 (p_1+2p_2)}=\frac{1}{p_1+2p_2}n p_2^2 \ge (1-p_2)n p_2^2.$$

This heuristic argument leads to the following lemma:

\begin{lemma}\label{lemmaE}
Let $p_1>1/2$ and let
$D_1\!\!:=\left\{LC_n\ge np_1 +\left((1-p_2)^2-p_2\right) n p_2^2 \right\}$. Then,
$\mathbb{P}(D_1) \ge 1- 4\exp(-2n p_2^6)-\exp \left(n(p_2^3+\log(1-p_2^3))(p_1-p_2^3)\right)$.
\end{lemma}
$\textbf{Proof.}$ For $p_1>\delta>0$, let $D_2^x(\delta) :=\left\{\left|\sum_{i=1}^n
\mathbf{1}_{\{X_i=\alpha_1\}}-np_1 \right|\le \delta n \right\}$, let
$D_2^y(\delta) :=\left\{\left|\sum_{i=1}^n \mathbf{1}_{\{Y_i=\alpha_1\}}-np_1 \right|
\le \delta n \right\}$,
and let $D_2(\delta):=D_2^x(\delta) \cap D_2^y (\delta)$, so that
on $D_2(\delta)$, at least $n_1(\delta):= n(p_1-\delta)$
letters $\alpha_1$ can be aligned.
Clearly, $1+\min(R_i^2, S_i^2)$ has a geometric distribution with
parameter $1-(p_2/{(p_1+p_2)})^2$.  Also, if $\mathcal{G}_1, \dots, \mathcal{G}_r$ are
iid geometric random variables with parameter $p$, then for any $\beta < 1$,
\begin{equation}\label{MatI}
\mathbb{P}\left(\sum_{i=1}^r \mathcal{G}_i \le \frac{\beta}{p}r\right)
\le \exp\left(-(\beta-1-\log\beta)r\right).
\end{equation}
By taking $p=1-(p_2/{(p_1+p_2)})^2$ and $r=n_1(\delta)$, and since the sequences have same 
length $n$, the following equality in law holds true:
$$\sum_{i=1}^{n_1(\delta)} \min(R_i^2, S_i^2) + n_1(\delta) \overset{d}{=}\sum_{i=1}^{n_1(\delta)} \left(\mathcal{G}_i \wedge n\right).$$
For any $\beta <1$, let us estimate
$$
\mathbb{P}\left(\sum_{i=1}^{n_1(\delta)} \min(R_i^2, S_i^2)
< \frac{\beta n_1(\delta)}{1-\left(\frac{p_2}{p_1+p_2}\right)^2}-
n_1(\delta) \right).
$$
First,
$$\frac{\beta n_1(\delta)}{1-\left(\frac{p_2}{p_1+p_2}\right)^2}-
n_1(\delta) \le n,$$
and therefore,
\begin{equation}\label{min}
\mathbb{P}\left(\sum_{i=1}^{n_1(\delta)} \min(R_i^2, S_i^2)
< \frac{\beta n_1(\delta)}{1-\left(\frac{p_2}{p_1+p_2}\right)^2}-
n_1(\delta) \right) \le e^{-(\beta-1-\log \beta)n_1(\delta)}.
\end{equation}

Next, let
$$D_3(\beta,\delta):=\left\{\sum_{i=1}^{n_1(\delta)} \min(R_i^2, S_i^2) \ge \frac{\beta n_1(\delta)}{1-\left(\frac{p_2}{p_1+p_2}\right)^2}-n_1(\delta)\right\}.$$
Letting $\delta=p_2^3$ and $\beta=1-p_2^3$, and 
when $D_2(\delta)$ and $D_3(\beta, \delta)$ both hold, then
\begin{align}
LC_n &\ge \frac{\beta n_1(\delta)}{1-\left(\frac{p_2}
{p_1+p_2}\right)^2}-n_1(\delta)+n_1(\delta)\notag\\
%    &=\frac{n_1(\delta)}{1-\left(\frac{p_2}{p_1+p_2}\right)^2}-n_1(\delta)+n_1(\delta)
%   -\frac{n_1(\delta) p_2^3}{1-\left(\frac{p_2}{p_1+p_2}\right)^2}\notag\\
     & =np_2^2\frac{p_1-p_2^3}{(p_1+p_2)^2-p_2^2}+n(p_1-p_2^3)-n p_2^2
     \frac{p_2 (p_1-p_2^3)}{1-\left(\frac{p_2}{p_1+p_2}\right)^2}\notag\\
     & =np_1 +\left(\frac{(p_1-p_2^3)(1-p_2(p_1+p_2)^2)}{p_1(p_1+2p_2)}-p_2\right) n p_2^2\notag\\
     & \ge np_1 +\left(\frac{(p_1-p_2^3)(1-p_2)}{p_1(1+p_2)}-p_2\right) n p_2^2\notag\\
     & \ge np_1 +\left((1-p_2)^2-p_2\right) n p_2^2.\notag
     \end{align}
Since $D_2(p_2^3)\cap D_3(1-p_2^3,p_2^3) \subset D_1$, it follows from 
Hoeffding's inequality and \eqref{min} that
$$\mathbb{P}(D_1) \ge 1- 4\exp(-2n p_2^6)-\exp \left(n(p_2^3+\log(1-p_2^3))(p_1-p_2^3)\right).$$
% \le 5 \exp\left(-\frac{ np_2^6}{4}\right).
 \CQFD

To state our next lemma, let us introduce some more notation. First, let
\begin{equation}\label{VK}
V(k):=\left\{(v_1, v_2, \dots, v_k) \in \mathbb{Z}^k: |v_1|+\dots +|v_k|\le 2k \right\},
\end{equation}
and then let
\begin{equation}
P:= \bigcup_{2 k\ge  n p_2^2 }  V(k).
\end{equation}
With these definitions, the previous lemma further yields:

\begin{lemma}\label{coro1}
Let $O$ be the set of all the optimal alignments of $X = (X_i)_{1\le i \le n}$
and $Y = (Y_i)_{1 \le i \le n}$, let $D=\{O \subset P \}$, let $p_1 > 1/2$
and let $p_2 < 1/10$. Then,
$\mathbb{P}(D)\ge 1-5 \exp\left(-{ np_2^6}/{5}\right)$.
\end{lemma}
$\textbf{Proof.}$ Let $N_1^{X}$ be the number of letters $\alpha_1$ in $X$,
and $N_1^Y$ be the corresponding number in $Y$, and so $N_1=N_1^X+N_1^Y$.
From the proof of the previous lemma, with its notation, it is clear that:
\begin{align}
D_1\cap D_2(p_2^3) &\subset \left\{LC_n \ge
\frac{N_1}{2}-n p_2^3 +\left((1-p_2)^2-p_2\right) n p_2^2\right\}\\
&\subset \left\{LC_n \ge \frac{N_1}{2}+\frac{1}{2} n p_2^2 \right\}:={\widetilde D_1}(p_2^2),
\end{align}
since $p_2<1/10$.  But, $D_2(p_2^3) \cap D_3(1-p_2^3, p_2^3) \subset D_1$, so as in the previous lemma,
\begin{align}
\mathbb{P}\!\left(\!LC_n \ge \frac{N_1}{2}+\frac{1}{2} n p_2^2\!\right)
&\ge 1- 4\exp(-2n p_2^6)-\exp \left(n(p_2^3+\log(1-p_2^3))(p_1-p_2^3)\right)\nonumber\\
&\ge 1-5 \exp\left(-{ np_2^6}/{5}\right), \nonumber
\end{align}
since again $p_2<1/10$.  It remains to show that ${\widetilde D_1}(p_2^2)\subset D$.
But, for any alignment with $|v|=k \ge 0$,
\begin{equation}\label{LCS:add0}
 LC_n\le \frac{N_1}{2}-\frac{1}{2}\sum_{i=1}^{k} |v_i|+ k,
\end{equation}
while on ${\widetilde D_1}(p_2^2)$,
\begin{equation}\label{LCS:add1}
LC_n \ge \frac{N_1}{2}+\frac{1}{2} n p_2^2.
\end{equation}
In case $|v|=0$, no optimal alignment do satisfy both \eqref{LCS:add0} and \eqref{LCS:add1},
while for $|v|\ge 1$, they both combine to yield
$\sum_{i=1}^k |v_i| \le 2k \text{ and } n p_2^2 \le 2 k$,
and this finishes the proof.
\CQFD

The previous lemma asserts that, with high probability, any optimal alignment belongs to
the set $P$.  Hence, in order to prove that the optimal alignments satisfy a property,
one needs, essentially, to only prove it for the alignments in $P$.

\subsection{High probability events}

Recall, from Definition~\ref{D:alignment}, that any
$v\in \mathbb{Z}^k, k\ge 1$ is associated with an alignment having
$k=|v|$ cells $C_v(1), \dots, C_v(|v|)$, and that a cell is called a
nonzero-cell if it contains a different number of
letters $\alpha_1$ on the $X$-strand and on the $Y$-strand.
For any $\theta > 0$, let $W^\theta$ be the subset of $P$, consisting
of the alignments having a proportion of nonzero-cells at least equal to $\theta$, {\it i.e.},
$$W^\theta:=\left\{v\in P : |\{i\in [1,k]: v_i \ne 0\}|\ge \theta |v|\right\},$$
and let $(W^\theta)^c:=P \backslash W^\theta$.

To complete the proof of the theorem, some further relevant events need to be defined.

\begin{itemize}
  \item
  For any $v \in P$, let $E_v^\theta$ be the event that the proportion of
  zero-cells in $C_v(1), \dots, C_v(|v|)$, is at least equal to $\theta$. Then, let
$$E^\theta:=\bigcap_{v \in (W^\theta)^c} E_v^\theta:=\bigcap_{v \in (W^\theta)^c} \left\{I_b\ge \theta J_0\right\},$$
where $J_0$ is the number of zero-cells while $I_b$ is the number of breakable
zero-cells for $v$,
{\it i.e.}, $E^\theta$ is the event that every $v \in (W^\theta)^c$ has a proportion of breakable
zero-cells at least equal to $\theta$.
  \item Recall also from \eqref{weakside} and \eqref{totalnon1}, that $N_v^-$ is the total number of
non-$\alpha_1$ letters in the cell strands with the lesser number of $\alpha_1$, and
that $N_{>1}$ is the total number of non-$\alpha_1$ letters in $X$ and $Y$.  Then, let
   $$F^\theta:=\bigcap_{v\in W^\theta}F_v:=\bigcap_{v\in W^\theta}\left\{N_v^-\ge \frac{K}{m} N_{>1}\right\},$$
  {\it i.e.}, $F^\theta$ is the event that for every $v\in W^\theta$, the
  proportion of non-$\alpha_1$ letters which
  are on the cell-strand with the smaller number of letters $\alpha_1$, is at least equal to $K/m$.
  \item Let
    $$G^\theta:=\bigcap_{v \in W^\theta} G_v:= \bigcap_{v \in W^\theta}\left\{2|v| \le \frac{K}{2m} N_{>1}\right\},$$
    {\it i.e.}, $G^\theta$ are the alignments $v \in W^\theta$ having a proportion of aligned non-$\alpha_1$ letters at most equal to $K /2m$.
\end{itemize}

Finally recall from Section~\ref{sectiondef} that $A_n=\{(X,Y)\in \mathcal{B}_n\}$ is the event
that there exists an optimal alignment, with $|v|\ge 1$,
such that $N_v^-\ge K N_{>1}/m$ and $2|v|\le K N_{>1}/2m$, and  therefore
\begin{equation}\label{finalrela}
D\cap E^\theta \cap F^\theta \cap G^\theta \subset A_n.
\end{equation}
Our next task is to prove that each one of the
events $E^\theta,F^\theta,G^\theta$ hold with high probability.
Let us start with $E^\theta$.

\begin{lemma}\label{lemma2}
Let $0<\theta \le p_1^2/(1+p_1^2)$, then
\begin{equation}\label{LCS:inLemma2}
\mathbb{P}(E^\theta) \ge 1- \sum_{2 k\ge  n p_2^2  }\exp\left(-\left(2(1-\theta) \left(\frac{p_1^2}{1+p_1^2}-\theta\right)^2-\log f(\theta)\right)k\right),
\end{equation}
where $f(\theta)=\left((4+2\theta)/{\theta^2}\right)^{\theta}\left((2+\theta)/{2}\right)^2
\left({1}/(1-\theta)\right)^{1-\theta}$.
\end{lemma}
$\textbf{Proof.}$ For any $v\in P \backslash W^\theta$, let us compute the probability that a zero-cell
in the alignment associated with $v$ is breakable.
Recalling the definition of $T_0$ in \eqref{stoppingtime}, for $2\le j \le m$, let $M_j$ be the event
that this cell ends with a pair of letters $\alpha_j$. So, when $M_j$ holds, then $T_0=T_0^j$.
For $2\le j \le m$, let also
$$U_1^j:=\min\{i=2,3,\ldots:\quad R_{i-1}^j\ne 0,\quad S_{i-1}^j=0,\quad R_{i}^j =0, \quad S_{i}^j \ne 0\},$$
$$U_2^j:=\min\{i=2,3,\ldots:\quad R_{i-1}^j= 0,\quad S_{i-1}^j\ne 0,\quad R_{i}^j \ne 0, \quad S_{i}^j = 0\},$$
and
$$ \quad U^j:=\min\{U_1^j,U_2^j\}.$$

With the above constructions, conditional on the event $M_j$, if $U^j <T_0^j $ then this
zero-cell is breakable and thus, to lower bound the probability that it is breakable,
it is enough to lower bound $\mathbb{P}(U^j <T_0^j )$.
To do so, let first $(Z_i^j)_{i\ge 1}$ be the independent random vectors given by:
$$Z_i^j=(R_{2i-1}^j,S_{2i-1}^j,R_{2i}^j,S_{2i}^j).$$
Then, let
$$\tilde{U}^{j}=\min\{i=1,2,\ldots: \quad Z_i^j \in B_1 \cup B_2\}, $$
$$\tilde{T}_0^j=\min \{i=1,2,\ldots: \quad Z_i^j \in B_3 \cup B_4\}, $$
where
$$B_1:=\mathbb{N}^* \times \{0\} \times \{0\}\times \mathbb{N}^* ,
B_2:=\{0\}\times \mathbb{N}^* \times \mathbb{N}^* \times \{0\},$$
$$B_3:= \mathbb{N}^* \times \mathbb{N}^* \times \mathbb{N} \times \mathbb{N}, \quad
B_4:=      \mathbb{N}\times \mathbb{N} \times \mathbb{N}^* \times \mathbb{N}^* ,$$
and where as usual $\mathbb{N}$ is the set of non-negative
integers, while $\mathbb{N}^*=\mathbb{N}\backslash \{0\}$.
Clearly, $2\tilde{U}^{j} \ge U^j$ and $2\tilde{T}_0^j -1 \le T_0^j$,
thus $\mathbb{P}(U^j <T_0^j) \ge \mathbb{P}(2\tilde{U}^{j}
< 2\tilde{T}_0^j -1)=\mathbb{P} (\tilde{U}^{j} < \tilde{T}_0^j)$.  Now, since the random vectors
$(Z_i^j)_{i\ge 1}$ are iid, and since $B_1\cup B_2$ and $B_3 \cup B_4$ are pairwise disjoint,
\begin{align}
\mathbb{P}(\tilde{U}^{j} < \tilde{T}_0^j)&=\frac{\mathbb{P}(Z_i^j\in B_1 \cup B_2)}
{\mathbb{P}(Z_i^j\in B_1 \cup B_2)+\mathbb{P}(Z_i^j\in B_3 \cup B_4)}\notag\\
&=\frac{2p_1 ^2}{2p_1 ^2+2(p_1+p_j)^2-p_j^2}
\ge  \frac{p_1^2}{1+p_1^2 }.\notag
\end{align}
Therefore,
\begin{align}
\mathbb{P}(\text{a zero-cell is breakable})&=\sum_{j=2}^m \mathbb{P}
(\text{a zero-cell is breakable}|M_j)\mathbb{P}(M_j)\notag\\
& = \sum_{j=2}^m \mathbb{P}(U^j <T_0^j)\mathbb{P}(M_j)
 \ge \frac{p_1^2}{1+p_1^2 }.\notag
\end{align}

Let $J$ be the index set of all the zero-cells in the alignment associated with $v\in (W^\theta)^c$,
and so $|J|\ge (1-\theta) |v|$.  For each $i \in J$, let $I_i$ be the Bernoulli
random variable which is one if the cell $C_v(i)$ is breakable and 0 otherwise.
Recall that $E_v^\theta$ is the event that the proportion of
breakable cells in $v$ is at least equal to $\theta$.
Then, since $\theta \le p_1^2/(1+ p_1^2)$, from Hoeffding's inequality, and after subtracting the mean,
\begin{align}
\mathbb{P}((E_v^\theta)^c)&=\mathbb{P}\left(\sum_{i \in J}I_i < \theta |J|\right)
%\notag\\
%&=\mathbb{P}\left(\sum_{i \in J}I_i - \mathbb{E}\left(\sum_{i \in J}I_i\right)
%< \theta |J|- \mathbb{E}\left(\sum_{i \in J}I_i\right)\right)\notag\\
%&
 \le \exp\left(-2(1-\theta) |v|\left(\frac{p_1^2}
{1+p_1^2}-\theta\right)^2\right).\notag
\end{align}
Recall now the definition of $V(k)$ in \eqref{VK} and let $(W^\theta(k))^c:=(W^\theta)^c \cap V(k)$.
For any two integers, $\ell$ and $q \ell $, with $0<q<1$, Stirling's formula in the form
$1 \le {\ell!e^\ell}/(\sqrt{2 \pi \ell} \ell^\ell) \le {e}/{\sqrt{2 \pi}}$,
gives
\begin{equation}\label{binomial}
\binom{\ell}{q \ell } \le q^{-q \ell} (1-q)^{-(\ell-q \ell)},
\end{equation}
which, when combined with simple estimates yields,
\begin{align}
|(W^\theta(k))^c| &\le 2^{\theta k}\binom {2k+\theta k}
{\theta k}\binom{k}{\theta k} \notag\\
&\le (f(\theta))^k 
:=\left(\left(\frac{4+2\theta}{\theta^2}\right)^{\theta}
\left(\frac{2+\theta}{2}\right)^2\left(\frac{1}{1-\theta}\right)^{1-\theta}\right)^k.
\end{align}
Next, let $E^\theta(k)=\bigcap_{v \in (W^\theta(k))^c} E_v^\theta$,
then 
\begin{align*}
\mathbb{P}((E^\theta(k))^c)&\le\sum_{v \in (W^\theta(k))^c}\mathbb{P}
((E^\theta_v)^c) \\ 
&\le \exp\!\left(-\left(\!2(1-\theta)\left(\frac{p_1^2}{1+p_1^2}-\theta\right)^2 
-\log f(\theta)\!\right)k\right),
\end{align*}
and therefore,
\begin{align}
\mathbb{P}((E^\theta)^c)&\le \sum_{2 k\ge  n p_2^2  }\mathbb{P}((E^\theta(k))^c) \notag\\
&\le \sum_{2 k\ge  n p_2^2  }\exp\left(-\left(2(1-\theta) \left(\frac{p_1^2}
{1+p_1^2}-\theta\right)^2-\log f(\theta)\right)k\right).\label{last}
\end{align}\CQFD

Of course, in \eqref{last}, one wants
\begin{equation}\label{LCS:in2}
2(1-\theta) \left(\frac{p_1^2}{1+p_1^2}-\theta\right)^2-\log f(\theta)>0,
\end{equation}
and choices of $\theta$ for which this is indeed the case are given later.

Let $u$ be a non-negative integer. For any $-u$-cell ending with an aligned pair of
letters $\alpha_j$ (the event $M_j$ holds for this cell),
let $\tau_X^j(\ell)$ be the index of the $\ell$-th $R_i^j$ such that $R_i^j\ne 0$, \emph{i.e.},
$$\tau_X^j(1)=\min\{i\ge 1 : \  R_i^j \ne 0\}, $$
and for any $\ell \ge 1$, $\tau_X^j(l+1)=\min \{i> \tau_X^j(\ell): \ R_i^j \ne 0 \}$.  Let
$$\rho^{j,-}:= \min\{\ell=1,2,\ldots: \  S_{u+\tau_X^j(\ell)}^j \ne 0\}.$$
In words, $\rho^{j,-}$ is the number of nonzero values taken by $R^j=(R_i^j)_{1 \le i \le s}$
(where $s$ is the number of letters $\alpha_1$ in the $X$-strand of the cell).  Since $X$ and $Y$ are independent,
\begin{align}
\mathbb{P}(\rho^{j,-}=k)&=\mathbb{P}(S_{u+\tau_X^j(1)}^j =
0,\dots, S_{u+\tau_X^j(k-1)}^j = 0, S_{u+\tau_X^j(k)}^j \ne 0)\notag\\
&=\left(\frac{p_1}{p_1+p_j}\right)^{k-1} \frac{p_j}{p_1+p_j},\label{Geometric}
\end{align}
for $k=1,2,\dots$.
Thus, $\rho^{j,-}$ has a geometric distribution with
parameter $\tilde{p}_j={p_j}/(p_1+p_j)$, $2\le j \le m$.
(By just replacing $\tau_X$ by $\tau_Y$ the random variables $\rho^{j,-}$ can then
be defined for $u$-cells.  Hence, since $X$ and $Y$ have the same law,
the corresponding law of $\rho^{j,-}$ remains unchanged,
therefore taking care of all the cases.)
When $-u<0$, the number of letters $\alpha_j$ in the X-strand
(which is the strand with the smaller number of letters $\alpha_1$) is
at least $\rho^{j,-}-1$ and, as shown in the next lemma,
this provides a lower bound for $N_v^-$ (the number of non-$\alpha_1$ letters on the
cell-strand with the lesser number of letters $\alpha_1$) in this $-u$-cell.

Recalling now that $F^\theta=\bigcap_{v\in W^\theta}\left\{N_v^-\ge K N_{>1}/m \right\}$, we have:
\begin{lemma}\label{lemma3}
Let $0<\theta<1$, let $K=2^{-4} 10^{-2} e^{-67}$, and let $p_1 \ge 1-  e^{-67}/4$. Then,
$\mathbb{P}(F^\theta)\ge 1-38 \exp(- 3n p_2^2/200)$.
\end{lemma}
$\textbf{Proof.}$ For any $v\in W^\theta$, let $J$ be the index set of all the nonzero-cells of the
alignment corresponding to $v$, hence, $|J|\ge \theta |v|$.  Then,
$$N_v^-=\sum_{i=1}^{|v|} N_v^-(i)= \sum_{i \in J} N_v^-(i)
\ge \sum_{i \in J} \left(\rho_i^{j(i),-}-1\right),$$
where $j(i)$ is the index of the last aligned pair of letters
$\alpha_j$ in the cell $C_v(i)$, and where $\rho_i^{j(i),-}$ is the number of
nonzero $R^{j(i)}=(R_{\ell}^{j(i)})_{1 \le \ell \le s}$
(assuming this is a $-u$-cell, and that $s$ is the number of letters $\alpha_1$ in the $X$-strand of $C_v(i)$. In case of a $u$-cell, by symmetry, the same argument is valid on the $Y$-strand).
From \eqref{Geometric}, $\rho_i^{j(i),-}$ is a geometric random variable
with parameter $\tilde{p}_{j(i)}$. Now, let $\varepsilon>0$,
let again $\tilde{p}_2={p_2}/(p_1+p_2)$, and let
$F_{1,v}:=\left\{N_v^-\ge \varepsilon|v|/\tilde{p}_2\right\}$.
Then,
\begin{align}
\mathbb{P}(F_{1,v}^c) &\le \mathbb{P}\left(\sum_{i \in J} \left(\rho_i^{j(i),-}-1\right)
\le \frac{\varepsilon}{\tilde{p}_2}|v| \right)\notag\\
 & \le \mathbb{P}\left(\sum_{i \in J} \rho_i^{j(i),-} \le
 \frac{{\varepsilon}/{\theta}+\tilde{p}_2}{\tilde{p}_2}|J| \right)\notag\\
 & \le \mathbb{P}\left(\sum_{i \in J} \rho_i^{j(i),-}
 \le \frac{\varepsilon/\theta+2p_2}{\tilde{p}_2}|J| \right).\label{estimate3}
\end{align}
The geometric random variables $\rho_i^{j(i),-}$, $i\in J$, are independent each with parameter
$\tilde{p}_{j(i)}\le \tilde{p}_2$, and moreover the sequences have finite length $n$, therefore,
$$\mathbb{P}\left(\sum_{i \in J} \rho_i^{j(i),-} \le \frac{\varepsilon/\theta+2p_2}
{\tilde{p}_2}|J| \right) \le \mathbb{P}\left(\sum_{i \in J}
(\mathcal{G}_i \wedge n) \le \frac{\varepsilon/\theta+2p_2}{\tilde{p}_2}|J| \right),$$
where the $\mathcal{G}_i$ are iid geometric random variables with
parameter $\tilde{p}_2$.  As proved later, and using \eqref{MatI}, when
\begin{equation}\label{exponentialPb1}
\frac{\varepsilon/\theta+2p_2}{\tilde{p}_2}|J| < n,
\end{equation}
it follows that
\begin{align}\notag
\mathbb{P}\left(\sum_{i \in J} \rho_i^{j(i),-}
\le \frac{\varepsilon/\theta+2p_2}{\tilde{p}_2}|J| \right)
&\le \mathbb{P}\left(\sum_{i \in J} \mathcal{G}_i
\le \frac{\varepsilon/\theta+2p_2}{\tilde{p}_2}|J| \right)\\
& \le \exp\left(\left(1 +\log ({\varepsilon}/{\theta}+2 p_2 )\right)
\theta |v|\right)\label{LCS:estimate1}.
\end{align}

\noindent
Let $F^\theta_{1}(k):=\bigcap_{v \in W^\theta \cap V(k)}F_{1,v}
=\bigcap_{v \in W^\theta \cap V(k)}\left\{N_v^-\ge \varepsilon|v|/\tilde{p}_2\right\}$, and let
$F^\theta_{1}:=\bigcap_{2 k\ge  n p_2^2} F^\theta_1(k)$.
From the very definition of $V(k)$ in \eqref{VK}, and using
\eqref{binomial},
$$|V(k)|\le 2^k \binom{3k}{k} \le 2^k 3^k \left(\frac{3}{2}\right)^{2k}=
\left(\frac{27}{2}\right)^k,$$
which when combined with \eqref{LCS:estimate1} leads to
\begin{align}
\mathbb{P}(F^\theta_{1}(k))&\ge
 1- \exp\left(k \log(27/2)+k\left(1 +\log ({\varepsilon}/{\theta}+2 p_2 )\right)\theta \right).
\end{align}
Of course, one wants $\log(27/2)+\left(1 +\log ({\varepsilon}/{\theta}+2 p_2 )\right)\theta  <0$.
Choosing $\theta=1/25$ and $\varepsilon = 10^{-2} e^{-67}$, then
$\mathbb{P}((F^\theta_{1}(k))^c)\le e^{-3k/100}$, for any $p_1 \ge 1- 2^{-2} e^{-67}$,
and so
$$\mathbb{P}((F^\theta_{1})^c)\le \sum_{2k \ge  n p_2^2} 
\mathbb{P}((F^\theta_{1}(k))^c)\le \frac{\exp(- 3n p_2^2/200)}{1- \exp(- 3/100)} \le 34
\exp(- 3n p_2^2/200).$$
Note also that for these choices of $\theta$ and $p_1$, \eqref{LCS:in2} is satisfied and
so $E^\theta$ also holds with high probability.

From the proof of Lemma~\ref{lemmaE}, when $D_2((1-p_1))$ holds, the total number of
non-$\alpha_1$ letters in $X$ and $Y$ is at most $4n(1-p_1)$.
Thus $N_{>1} \le 4 n (1-p_1)$, and so
when $F^\theta_1 \cap D_2((1-p_1))$ holds, for every $v\in W^\theta$,
$$\frac{N_v^-}{N_{>1}}\ge \frac{\varepsilon |v| }{\tilde{p}_2 4n (1-p_1)}
\ge \frac{\varepsilon}{\tilde{p}_2 4n (1-p_1)} \frac{ n p_2^2}{2}
\ge \frac{\varepsilon p_2}{16 (1-p_1)}\ge \frac{\varepsilon}{16 m} \ge \frac{K}{m}.$$
Also note that by properly choosing these constants and under the further
condition $400 m K < 1$, it follows that
\eqref{exponentialPb1} holds true.  Therefore,
\begin{align}\notag
\mathbb{P}((F^\theta)^c)&\le \mathbb{P}((F^\theta_1)^c)+\mathbb{P}( (D_2(1-p_1))^c)\\
&\le 34
\exp(- 3n p_2^2/200) + 4\exp(-2n(1-p_1)^2)\notag \\
&\le 38 \exp(- 3n p_2^2/200)\notag.
\end{align}\CQFD

Recalling that $G^\theta = \bigcap_{v \in W^\theta}\left\{2|v| \le K N_{>1}/{2m}\right\}$, 
we finally have:
\begin{lemma}\label{lemma4}
Let $0<\theta<1$, let $K=2^{-4} 10^{-2} e^{-67}$, and moreover let $p_2 \le \min\{2^{-2}e^{-5}\!K/m, K/2m^2\}$.  Then, $\mathbb{P}(G^\theta)\ge 1-8 \exp(- n p_2^2/2)$.
\end{lemma}
$\textbf{Proof.}$  For any $v \in W^\theta $, let $C_v(1), \dots, C_v(|v|)$ be the corresponding cells.
If the cell $C_v(i)$ ends with a pair of aligned $\alpha_j$, $2\le j \le m$, then
let $\rho_i^{j(i)}$ be the number of nonzero values
taken by $R^{j(i)}$ in $C_v(i)$.  If $v_i\le 0$, by the same arguments as in
getting \eqref{Geometric}, $\rho_i^{j(i)}$ has a geometric distribution
with parameter $\tilde{p}_{j(i)} = p_{j(i)}/(p_1 + p_{j(i)})$.  If $v_i >0$, then 
there exists a geometric
random variable $\rho_i^{j(i),-}$ with parameter $\tilde{p}_{j(i)}$ such
that $\rho_i^{j(i),-}\le \rho_i^{j(i)}\le \rho_i^{j(i),-}+v_i$.
Let $N_{>1}^X$ (resp. $N_{>1}^Y$) be the number of non-$\alpha_1$ letters in $X$
(resp. $Y$), so that $N_{>1}=N_{>1}^X+N_{>1}^Y$, and let
$$G_v^X:=\left\{|v| \le \frac{K}{2m} N_{>1}^x\right\}\text{ and }G_v^Y:=\left\{|v| \le \frac{K}{2m}
N_{>1}^y\right\},$$
and so $G_v^X \cap G_v^Y \subset G_v$.  Since $N_{>1}^X \ge \sum_{i=1}^{|v|} \rho_i^{j(i)}$,
\begin{align}
\mathbb{P}\left((G_v^X)^c\right)& \le  \mathbb{P}\left(|v|>
\frac{K}{2m}\sum_{i=1}^{|v|} \rho_i^{j(i)} \right)\notag\\
&\le  \mathbb{P}\left(|v|> \frac{K}{2m}\left(\sum_{1\le i\le |v|, v_i \le 0}
\rho_i^{j(i)}+\sum_{1\le i\le |v|, v_i > 0} \rho_i^{j(i),-} \right)\right)\notag\\
& \le \mathbb{P}\left(  \sum_{i=1}^{|v|} (\mathcal{G}_i \wedge n)<\frac{2m|v|}{K}\right ),\notag
\end{align}
where the $\mathcal{G}_i$ are iid geometric random variables with parameter
$\tilde{p}_2$ and the truncation is at $n$ since the sequences have such a length. 
From the proof of Lemma~\ref{lemmaE}, when $D_2((1-p_1))$ holds, 
$N_{>1} \le 4 n (1-p_1)$, then $|v| \le 2 n(1-p_1)$.  
Thus $2m|v| \le 2mn(1-p_1) < 2m^2p_2n$, and so if
${2m^2 p_2}<K$, then for any $p_2\le 2^{-2} e^{-5} K/m$,

\begin{align}
&\mathbb{P}\left((G_v^X)^c \cap D_2((1-p_1))\right) \notag\\
 &\le \mathbb{P}\left(  \sum_{i=1}^{|v|} \mathcal{G}_i <\frac{2m|v|}{K}\right )
 \le \mathbb{P}\left( \sum_{i=1}^{|v|} \mathcal{G}_i  <
\frac{e^{-5}|v|}{\tilde{p}_2} \right)
 \le \exp(-4|v|).
\end{align}
Likewise, $\mathbb{P}\left((G_v^Y)^c \cap D_2((1-p_1))\right) \le \exp(-4|v|),$ and
thus
$$\mathbb{P}\left((G_v)^c \cap D_2((1-p_1))\right) \le 2 \exp(-4|v|).$$
As before, let
$G^\theta(k):=\bigcap_{v \in W^\theta \cap V(k)} G_v \text{ and } 
G^\theta=\bigcap_{2 k\ge  n p_2^2} G^\theta(k)$,
then
$$\mathbb{P}((G^\theta(k))^c \cap D_2((1-p_1)) )\le |V(k)|2 \exp(-4k) \le 2 \exp(-k),$$
and
\begin{align}
&\mathbb{P}((G^\theta)^c) \le \mathbb{P}((G^\theta)^c \cap D_2((1-p_1))) 
+ \mathbb{P}(D_2((1-p_1))^c) \notag\\
&\le \sum_{2k \ge  n p_2^2 }\mathbb{P}((G^\theta(k))^c \cap D_2((1-p_1))) + 4\exp(-2n(1-p_1)^2)\notag\\
&\le \frac{2}{1-1/e} \exp(- n p_2^2/2) + 4\exp(-2n(1-p_1)^2) \notag\\
&\le 8 \exp(- n p_2^2/2).
\end{align}
\CQFD

From Lemma~\ref{coro1}--\ref{lemma4},
using \eqref{finalrela}, letting
$\theta=1/25$, $K= 2^{-4}10^{-2} e^{-67}$ and $K_m:=\min(K, 1/800m)$, and 
for
 $p_2 \!\le \! \min\{2^{-2}e^{-5}\!K_m/m, K_m/2m^2\}$, it follows that:  

\begin{align}
\mathbb{P}(A_n^c)&\le \mathbb{P}(D^c)+\mathbb{P}((E^\theta)^c)
+\mathbb{P}((F^\theta)^c)+\mathbb{P}((G^\theta)^c)\notag\\
&\le 5\exp\!\left(\!-\frac{ np_2^6}{5}\right)+74\exp\!\left(\!-\frac{n p_2^2}{10^3}\right)
+ 38\exp\!\left(\!-\frac{3n p_2^2}{200}\!\right)+8\exp\!\left(\!- \frac{n p_2^2}{2}\right)\notag\\
&\le 125\exp\!\left(\!-\frac{ np_2^6}{5}\right).
\end{align}
This finishes the proof of Theorem~\ref{T:theorem2}.\CQFD
\begin{remark}
(i) Our results on the central $r$-th absolute moments of the LCS continue to be valid for
three or more sequences of random words. First, the upper bound methods are 
very easily adapted to provide the
same order $n^{r/2}$.  Next, for the lower bound, the alignments can still be represented
with a series of cells, each of the cells ending with the same non-$\alpha_1$ letter
from every strand.  Then, with exponential bounds techniques, a similar high probability event can be exhibited,
also leading to a lower bound of order $n^{r/2}$.

(ii) With the methodology developed here, the results of \cite{AHM} and \cite{HM} can
also be generalized,
beyond the variance or the Bernoulli case, to centered absolute moments, $m$-letters alphabets and even to a general scoring function framework with scoring functions satisfying bounded differences conditions.
\end{remark}

% ------------------------------------------------------------------------

\subsection*{Acknowledgment}
Many thanks to Ruoting Gong and an anonymous referee for their detailed reading and numerous comments on this paper.

% ------------------------------------------------------------------------
\end{document}